\documentclass{article}
\usepackage[utf8]{inputenc}

\usepackage[english]{babel}

\usepackage[letterpaper,top=2cm,bottom=2cm,left=3cm,right=3cm,marginparwidth=1.75cm]{geometry}

\usepackage{mathtools}
\usepackage{changepage}
\usepackage{amsmath, amsfonts,amssymb}
\usepackage{multirow, multicol, makecell, booktabs, graphicx}
\usepackage{subcaption, natbib, authblk, footnote, url}
\usepackage[acronym,nogroupskip,nopostdot]{glossaries}
\usepackage{xltabular}
\usepackage[colorlinks=true, allcolors=blue]{hyperref}
\usepackage[short]{optidef}
\usepackage{subfiles, nameref, hyperref}
\hypersetup{
     colorlinks = true,
     citecolor  = blue,
     linkcolor = blue,
     urlcolor = blue
}
\usepackage{cleveref}
\usepackage{geometry}
\usepackage[linesnumbered,ruled,vlined]{algorithm2e}
\RestyleAlgo{ruled}
\usepackage{setspace}
\usepackage{bm}
\usepackage{xcolor}

\usepackage[title]{appendix}
\usepackage{array}

\newcommand{\pb}{{\sc JRTPCSSP}}
\newcommand{\algo}{{\sc AGGNNI-CG}}

\newcommand{\revision}[1]{{\color{black}#1}}
\newcommand{\revisiontwo}[1]{{\color{black}#1}}
\newcolumntype{L}[1]{>{\raggedright\let\newline\\\arraybackslash\hspace{0pt}}m{#1}}
\newcolumntype{C}[1]{>{\centering\let\newline\\\arraybackslash\hspace{0pt}}m{#1}}

\title{Boosting Column Generation with Graph Neural Networks \\
for Joint Rider Trip Planning and Crew Shift Scheduling}

\author{Jiawei Lu}
\author{Tinghan Ye\thanks{Corresponding author: \href{mailto:joe.ye@gatech.edu}{joe.ye@gatech.edu}}
}
\author{Wenbo Chen}
\author{Pascal Van Hentenryck}

\affil{H. Milton Stewart School of Industrial \& Systems Engineering,\protect\\ Georgia Institute of Technology, Atlanta, GA}
\date{}

\newacronym{mip}{MIP}{Mixed Integer Program}
\newacronym{ilp}{ILP}{Integer Linear Program}
\newacronym{drt}{DRT}{Demand-Responsive Transport}
\newacronym{maas}{MaaS}{Mobility-as-a-Service}
\newacronym{brt}{BRT}{Bus Rapid Transit}
\newacronym{cat}{CAT}{Chatham Area Transit}
\newacronym{nsf}{NSF}{National Science Foundation}

\begin{document}
\maketitle

\begin{abstract}

Optimizing service schedules is pivotal to the reliable, efficient,
and inclusive on-demand mobility. This pressing challenge is further
exacerbated by the increasing needs of an aging population, the
oversubscription of existing services, and the lack of effective
solution methods.  This study addresses the intricacies of service
scheduling, by jointly optimizing rider trip planning and crew
scheduling for a complex dynamic mobility service. The resulting
optimization problems are extremely challenging computationally for
state-of-the-art methods.

To address this fundamental gap, this paper introduces the Joint Rider
Trip Planning and Crew Shift Scheduling Problem (\pb{}) and a novel
solution method, called Attention and Gated GNN-Informed 
Column Generation (\algo{}), that hybridizes column generation and
machine learning to obtain near-optimal solutions to the \pb{} with
real-life constraints of the application. The key idea of the
machine-learning component is to dramatically reduce the number of
paths to explore in the pricing problem, accelerating the most
time-consuming component of the column generation. The machine
learning component is a graph neural network with an
attention mechanism and a gated architecture, which is particularly
suited to cater for the different input sizes coming from daily
operations.

\algo{} has been applied to a challenging, real-world dataset from the
Paratransit system of Chatham County in Georgia. It produces
substantial improvements compared to the baseline column generation
approach, which typically cannot produce high-quality feasible
solutions in reasonable time on large-scale complex instances.
\algo{} also produces significant improvements in service quality
compared to the existing system.


\end{abstract}

\clearpage
\section{Introduction}
\label{sec:intro}

Improving the efficiency, accessibility, and reliability of 
transportation service scheduling is a vital challenge in both 
academic research and practical applications. The increasing 
complexities of urban mobility, driven by growing traffic 
congestion, the environmental impact of private vehicles, 
and the operational rigidity of traditional transit systems, 
demand innovative solutions. Additionally, as the population 
ages and the demand for accessible transportation services 
rises, there is an urgent need to optimize how mobility 
services are scheduled and delivered.

Transportation service scheduling plays a pivotal role in 
ensuring that diverse rider needs are met efficiently. From 
managing public transit networks to coordinating on-demand 
services, effective scheduling is essential for improving 
operational efficiency and user satisfaction. The challenge 
becomes even more significant when accounting for dynamic 
real-world conditions such as fluctuating demand, variable 
travel times, and resource constraints, all of which can 
hinder the performance of traditional scheduling methods.

As another example, Paratransit services, designed to serve riders
with disabilities or those unable to access regular transit, provide a
compelling case for the complexities involved in transportation
service scheduling. These systems often experience oversubscription
due to the increasing needs of an aging population, coupled with
limited operational resources. This situation underscores the
importance of optimizing schedules to balance service quality with
operational costs.

\subsection{Motivations and Contributions}

The scheduling in these diverse transportation systems --- whether in 
traditional transit, on-demand services, or specialized Paratransit 
operations --- is a multifaceted challenge. Despite notable progress 
in service scheduling, several key challenges still require 
significantly more attention.
First, there is a discernible need for an integrated framework that 
simultaneously addresses {\em rider trip plans} and {\em crew shift 
scheduling}. To date, within the academic sphere and in practice, 
these elements are predominantly examined in isolation, highlighting 
the necessity for a more holistic approach.
\revision{Second, the inherent complexity of such an integrated problem in 
transportation systems poses significant computational challenges. 
This complexity arises from the need to jointly optimize intricate 
rider trip itineraries (often involving multiple requests per rider) 
alongside dynamically determining optimal crew shifts within maximum 
working hour limits, unlike many traditional routing problems with 
predefined schedules.
Third, evaluating solution methods often relies on synthesized datasets, 
leaving a void in capturing the realities encountered in the field. 
Particularly, these realities include complex and fluctuating demand 
patterns --- regarding spatial-temporal distributions and correlations 
between a rider's multiple requests --- factors that synthetic data 
often struggle to fully replicate, thereby impacting algorithm performance 
evaluation.}

This study aims at addressing these challenges by introducing an
integrated approach for the Joint Rider Trip Planning and Crew Shift
Scheduling Problem (\pb{}) for mobility systems. Moreover, to address the
computational challenges of the \pb{}, the study combines a
traditional column generation with an Attention-Based, Gated Graph
Neural Network to reduce the search space and find high-quality
solutions quickly. More precisely, the primary contributions of this
research can be summarized as follows:

\revision{
\begin{itemize}
\item Problem Formulation: This study formally defines and characterizes 
    the \pb{}, explicitly capturing the joint optimization of rider trip 
    planning and crew shift scheduling with relevant real-world constraints, 
    highlighting its complexities compared to related but distinct routing 
    and scheduling problems.

\item Efficient Solution Method: This study proposes \revisiontwo{Attention and Gated GNN-Informed 
Column Generation (\algo{})}, a method that 
    synergistically integrates column generation with a purpose-built Graph Neural Network (GNN). 
    The core idea is to leverage the GNN to effectively learn complex patterns 
    from the \pb{} pricing problem's structure and constraints, thereby 
    drastically reducing the search space by predicting the most promising 
    paths (columns) to explore.

\item Tailored GNN Architecture: This study designs a specific GNN 
    architecture incorporating attention mechanisms and gated units. This 
    architecture is motivated by the need to effectively process the 
    graph-structured pricing problem inherent in \pb{} and handle varying 
    input sizes typical of daily operations. Unlike existing approaches that 
    rely on local edge-level predictions or maintain computationally expensive 
    edge embeddings, the proposed model integrates edge features directly 
    into the attention weight calculations, eliminating memory bottlenecks 
    while preserving representational power. This architectural innovation 
    enables application to large-scale real-world instances and enhances 
    edge classification performance.

\item Validation on Real-World, Large-Scale Instances: This study evaluates
    the proposed \algo{} method through extensive computational experiments 
    using large-scale instances derived directly from real-world operational 
    data. This evaluation strategy is crucial, as the performance of machine
    learning models can be sensitive to the characteristics of the input data; 
    results on synthetic instances may not fully reflect performance on 
    operational challenges where demand patterns are complex and not easily 
    replicated.
\end{itemize}
}

\subsection{Related Work}

\subsubsection{Service Scheduling in Transportation Systems}

Service scheduling in transportation systems involves coordinating 
both the rider and crew aspects. This includes rider trip plan 
scheduling and crew shift scheduling, both of which have been 
widely studied in the literature due to their critical role in 
real-world operations.

For rider trip plan scheduling, the objective is to generate an
optimal travel plan that fulfill the travel needs of a rider from
origin to destination. A travel plan may consist of multiple travel
segments, each of which is specified by travel mode, departure time,
arrival time, and other information. Of course, rider trip requests
are fulfilled by drivers/vehicles, especially in on-demand
transportation systems. Therefore, trip plan scheduling is highly
correlated with driver/vehicle tasks, i.e., at what time to serve
which riders in what order. In the literature, the type of problem is
typically modeled as Dial-A-Ride
Problem (DARP), which is a variant of
the general vehicle routing problem. Readers are referred to for a
comprehensive review of the DARP by \cite{cordeau2007dial} and
\cite{ho2018survey}. Solution methods in the existing studies can be
classified into two classes: exact methods such as branch-and-cut
\citep{cordeau2006branch, ropke2007models}, branch-and-price
\citep{garaix2010vehicle,garaix2011optimization}, and
branch-and-price-and-cut algorithm \citep{qu2015branch,
  gschwind2015effective}; heuristics such as tabu search
\citep{cordeau2003tabu,kirchler2013granular}, simulated annealing
\citep{braekers2014exact}, large neighborhood search
\citep{ropke2006adaptive,DBLP:conf/cp/JainH11,gschwind2019adaptive},
and genetic algorithm
\citep{jorgensen2007solving,cubillos2009application}. In general,
exact methods take longer time to produce high-quality solutions with
bound information, while heuristics usually generate relatively good
solutions without bound information within a shorter time. In the
literature, most studies use either exact methods or heuristics; few
study combines exact methods and heuristics in solving DARP.

For crew shift scheduling in public transit systems, a key challenge
is the driver shift scheduling, whose objective is to design optimal
driver shifts based on driver availability and time-varying rider
travel needs in such a way that supply and demand are matched over
time. The literature on driver scheduling in public transportation
systems primarily focuses on traditional fixed-route systems (e.g.,
buses) and several key areas of research and methodological
developments. Notable studies include exploring genetic algorithms for
shift construction \citep{wren1995genetic}, defining efficient driver
scheduling methodologies \citep{toth2013efficient}, and integrating
vehicle and crew scheduling with driver reliability
\citep{andrade2021vehicle}. Further research delves into robust and
cost-efficient resource allocation for vehicle and crew scheduling,
addressing operational disruptions \citep{amberg2019robust}, and
developing new mathematical models for the Drivers Scheduling Problem
(DSP) that accurately reflect real-world complexities
\citep{portugal2009driver}. Studies also extend to addressing
scheduling problems with mealtime windows and employing Integer Linear
Programming (ILP) models for effective scheduling \citep{kang2019bus},
and integrating duty scheduling and rostering to enhance driver
satisfaction in public transit \citep{borndorfer2017integration}.
However, there is a notable gap in research focusing on emerging
on-demand transit systems, as most studies primarily concentrate on
conventional fixed-route bus systems, leaving room for exploration in
the context of modern, demand-responsive transport models.

\subsubsection{Machine Learning for Combinatorial Optimization}

The integration of machine learning techniques to address
combinatorial optimization problems has become an increasingly
prominent field of study in recent years. Interested readers can refer
to comprehensive reviews of recent progress by
\cite{bengio2021machine}, \cite{kotary2021end}, \cite{mazyavkina2021reinforcement}, and
\cite{karimi2022machine}. Research endeavors in this domain
generally fall into one of the two categories: the application of
standalone machine learning models to derive solutions for
combinatorial problems and the enhancement of traditional mathematical
optimization methods through machine learning.

In the first category, deep learning-based neural network 
architectures are typically designed to identify patterns 
that connect problem instances to their optimal solutions. 
This process may involve learning from a dataset of
existing solutions --- a supervised learning approach (e.g.,
\cite{vinyals2015pointer}) --- or discovering strategies through a
process of trial and error, akin to reinforcement learning (e.g.,
\cite{nazari2018reinforcement}) or both as in
\cite{DBLP:journals/jair/YuanCH22}. The goal is to train machine
learning models in a sufficiently robust manner to generate solutions
autonomously. Despite the progress, achieving high-quality solutions
solely with machine learning remains a challenging
frontier. Combinatorial optimization problems are inherently complex
and encompass vast solution spaces that are difficult to navigate
efficiently. As a result, while machine learning models have shown
promise, they often struggle to match the solution quality of
established optimization methods. Recent studies underscore the
persistent difficulty in bridging this quality gap, suggesting the
need for more sophisticated models and training techniques.

In the latter category, machine learning models, once adequately
trained, are often employed to streamline or assist with the most
time-intensive aspects of traditional mathematical optimization
methods. This can involve pivotal tasks such as choosing cutting
planes in branch-and-cut algorithms or selecting columns in
branch-and-price algorithms. By synergizing the strengths of both
mathematical optimization and machine learning, this hybrid approach
has demonstrated considerable promise in tackling a range of complex
problems. Specifically, in the realm of enhancing the column
generation algorithm through machine learning techniques - a critical
component in vehicle routing and service scheduling within transportation
systems - researchers have recorded notable
advancements. \cite{morabit2021machine} applied a learned model to
select promising columns from those generated at each iteration of
column generation to reduce the computing time of reoptimizing the
restricted master problem. \cite{shen2022enhancing} designed an ML
model to predict the optimal solution of pricing subproblems, which is
then used to guide a sampling method to efficiently generate
high-quality columns.

Further advances in this field involve using
machine learning to streamline pricing subproblems. Pioneering studies
have employed machine learning techniques to simplify underlying
graphs, thereby lessening the complexity of these
subproblems. \cite{morabit2023machine} utilized a random forest model
to predict edges likely to be part of the master problem's
solutions. Owing to the complexity and variable sizes of the graphs in
subproblems, their model bases its predictions on local rather than
global graph information. In addressing varying graph sizes,
\cite{yuan2022neural} developed a graph neural network model with
residual gated graph convolutional layers, inspired by
\cite{joshi2019efficient}. Their model evaluates the likelihood of
each edge being part of the final solution, leading to the
construction of a reduced graph that retains edges with higher
predicted probabilities for pricing subproblems. The method
proposed in this paper differs from \cite{yuan2022neural} in three key
ways: First, this study tackle both crew shift scheduling and rider trip
planning, unlike \cite{yuan2022neural} which focused only on the
former. Second, the proposed method is applied to considerably larger and more
complex real-world scenarios, whereas \cite{yuan2022neural} tested
their algorithm on instances of smaller scales. Third, contrary to
\cite{yuan2022neural}, the proposed architecture in this paper avoids
updating the edge embeddings in each layer, a process that is computationally 
inefficient and leads to prohibitive memory consumption at larger scales during training. Instead, the proposed model integrates
graph attention and residual gating mechanisms, incorporating
edge-level features in the attention weight calculations without the
need to keep track of the edge embeddings. This reduces the risk of memory overload and improves scalability.

\subsection{Outline of the Paper}

The rest of this paper is organized as
follows. Section~\ref{sec:problem} defines the \pb{}.
Section~\ref{sec:column_generation} presents the column generation
algorithm. Section~\ref{sec:gnn} describes the machine learning
methodology and its integration within the column generation method.
Section~\ref{sec:experiment} reports the benefits of the proposed
approach on the real-world dataset. Section~\ref{sec:conclusion}
concludes the paper and discusses future research directions.

\section{Joint Rider Trip Planning and Crew Shift Scheduling}
\label{sec:problem}

\subsection{Problem Description}

The \pb{} considered in this paper is a variant of the Dial-A-Ride
Problem (DARP). Compared to the classic DARP, the \pb{} introduces two
additional features. First, it considers scenarios where a single
rider may have multiple trip requests in a day and imposes a complete
service constraint: either all of the rider requests are fulfilled, or
none at all. This approach, which disallows partial servicing, is
motivated by the need to guarantee return trips (e.g., for dialysis
patients).  Second, the \pb{} enforces a maximum working duration for
each driver (e.g., 8 hours), but does not predefine driver shifts.
\revision{The \pb{} assumes a one-to-one correspondence between drivers and vehicles, which reflects realistic scenarios in practice. For ease of notation, the terms ``driver" and vehicle" are used interchangeably throughout the paper.}

Formally, the problem is defined as follows. Consider a set of riders
$U$, with each rider $u \in U$ associated with a set of trip requests
$R_u$. These requests are to be accommodated by a homogeneous fleet of
vehicles, denoted by $F$. Each trip request $r$ in $R_u$ has an
origin, a destination, and strict time windows for both departure and
arrival times. Although fulfilling every trip request is not
mandatory, it is crucial to note that partial servicing of a rider
requests is not allowed: either all requests in $R_u$ are served or
none of them. The complete set of trip requests from all riders is
denoted as $R = \bigcup_{u \in U} R_u$. Each vehicle $f$ in the fleet
$F$ has a capacity of $C$, and the fleet size is constant at $\lvert F
\rvert$. The working shifts of these vehicles are not predefined;
instead, they need to be strategically determined alongside the
scheduling of the trip requests. 
The earliest shift start time and latest shift end time are $\xi_s$
and $\xi_e$, respectively. The maximum number of working hours
of each vehicle is denoted by $\eta$. All vehicles begin and end their
service at a common depot. The travel times between the depot,
origins, and destinations are known and constant. The primary goal is
to design an optimal schedule for the vehicle working hours and the
service of trip requests that maximizes the number of requests
served. This objective is motivated by a real system in the field,
where the number of requests often exceeds the capacity of the
service.

\revision{Notations used in this section are provided in Table~\ref{tab:notation1}.}

\begin{table}[!ht]
    \centering
    \caption{\revision{Summary of Notation}}
    \label{tab:notation1}
    \begin{tabular}{L{1.8cm} L{13cm}}
        \hline
        \textbf{Symbol} & \textbf{Description} \\
        \hline
        \multicolumn{2}{l}{\textbf{Sets}} \\
        $U$ & Set of riders \\
        $R_u$ & Set of trip requests for rider $u \in U$ \\
        $R$ & Set of all trip requests, $R = \bigcup_{u\in U} R_u$ \\
        $F$ & Set of homogeneous vehicles \\
        $P/D$ & Set of pickup/drop-off nodes, one for each request $r \in R$ \\
        $P_u$ & Set of pickup nodes associated with rider $u \in U$ \\
        $N$ & Set of nodes in the graph $g$ for the arc-based model; $N = \{0\} \cup P \cup D \cup \{2n+1\}$ \\
        $E$ & Set of edges in the graph $g=(N,E)$ \\
        $\Omega$ & Set of feasible vehicle routes (paths) \\

        \multicolumn{2}{l}{\textbf{Parameters}} \\
        $n$ & Total number of trip requests; $n = |R|$ \\
        $C$ & Capacity of each vehicle \\
        $\xi_s/\xi_e$ & The earliest shift start time / latest shift end time \\
        $\eta$ & Maximum working duration (in hours) for each vehicle/driver \\
        $d_i$ & Demand at node $i \in N$; $d_i > 0$ for pickup, $d_i < 0$ for drop-off, $d_0=d_{2n+1}=0$ \\
        $s_i$ & Service time required at node $i \in N$ \\
        $[a_i, b_i]$ & Time window (earliest/latest service start time) for node $i \in N$ \\
        $t_{ij}$ & Travel time between node $i$ and node $j$ for $(i,j) \in E$ \\
        $\alpha_{r\theta}$ & Binary parameter: 1 if request $r$ is served by route $\theta \in \Omega$, 0 otherwise \\

        \multicolumn{2}{l}{\textbf{Decision Variables}} \\
        $x_{ij}^f$ & Binary variable: 1 if vehicle $f \in F$ traverses edge $(i,j) \in E$, 0 otherwise \\
        $z_u$ & Binary variable: 1 if rider $u \in U$ is served (all requests in $R_u$ fulfilled), 0 otherwise \\
        $T_i^f$ & Continuous variable: arrival time of vehicle $f \in F$ at node $i \in N$ \\
        $Q_i^f$ & Continuous variable: load in vehicle $f \in F$ upon departure at node $i \in N$ \\
        $y_r$ & Binary variable: 1 if trip request $r \in R$ is served, 0 otherwise \\
        $\lambda_{\theta}$ & Binary variable: 1 if route $\theta \in \Omega$ is selected, 0 otherwise \\
        \hline
    \end{tabular}
\end{table}

\subsection{An Arc-Based Model}
\label{sec:arc_model}

Figure~\ref{fig:arc_model} presents an arc-based model for the
\pb{}. The model is defined on the graph $g=(N,E)$, where $N$ denotes
the set of nodes and $E$ the set of edges. Each trip request $r \in R$
is represented by a pickup node $i$ and a corresponding drop-off node
$n+i$, included in the pickup node set $P$ and the drop-off node set $D$,
respectively, where $n = |R|$. Additionally, an origin depot node $0$
and a destination depot node $2n+1$ are created for the physical
depot, leading to the definition $N = \{0\} \cup P \cup D \cup
\{2n+1\}$. The set $E$ comprises edges that connect nodes in $N$, with
all connections subject to time window constraints. Each node $i \in
N$ is associated with a demand $d_i$, a service time $s_i$, and a time
window $[a_i, b_i]$ that corresponds to the earliest and latest
service start time at node $i$. \revision{Note that $a_0 = a_{2n+1}$ and $b_0 =
b_{2n+1}$ represent the earliest shift start time $\xi_s$ and the latest
shift end time $\xi_e$, respectively.} For each edge $(i,
j) \in E$, the travel time between nodes $i$ and $j$ is denoted by
$t_{ij}$.

\begin{figure}[!ht]
    \begin{maxi!}
%
	{}
%
	{\sum_{f \in F} \sum_{i \in P} x^f(\delta^+(i)) \label{eq:arc_model_obj}}
%
	{\label{formulation:arc_model}}
%
	{}
%
	\addConstraint
	{\sum_{f \in F} x^f(\delta^+(i))}
	{= z_u \quad \label{eq:request}}
	{\forall u \in U, \forall i \in P_u}
	\addConstraint
	{x^f(\delta^+(i)) - x^f(\delta^+(n+i))}
	{= 0 \quad \label{eq:same_vehicle}}
	{\forall i \in P, \forall f \in F}
	\addConstraint
	{x^f(\delta^+(i)) - x^f(\delta^-(i))}
        {= \begin{cases} 1 &\textrm{if } i = 0 \\ -1 &\textrm{if } i=2n+1 \\ 0 & i \in P \cup D \end{cases} \quad \label{eq:vehicle_flow_balance}}
	{\forall i \in N, \forall f \in F}
	\addConstraint
	{(T^f_i + s_i + t_{ij})x^f_{ij}}
	{\leq T^f_j \quad \label{eq:arrival_time}}
	{\forall (i,j) \in E, \forall f \in F}
	\addConstraint
	{a_i \leq T^f_i}
	{\leq b_i \quad \label{eq:time_window}}
	{\forall i \in N, \forall f \in F} 
	\addConstraint
	{T_{2n+1}^f - T_0^f}
	{\leq \eta \quad \label{eq:shift}}
	{\forall f \in F} 
	\addConstraint
	{(Q^f_i + d_j)x^f_{ij}}
	{\leq Q^f_j \quad \label{eq:onbard_passengers}}
	{\forall (i,j) \in E, \forall f \in F}
	\addConstraint
	{\max(0,d_i) \leq Q^f_i}
	{\leq \min(C,C+d_i) \quad \label{eq:capacity_bound_ab}}
	{\forall i \in N, \forall f \in F}
	\addConstraint
	{x^f_{ij}}
	{\in \{0,1\} \quad \label{eq:arc_model_x}}
	{\forall (i,j) \in E, \forall f \in F}
	\addConstraint
	{z_u}
	{\in \{0,1\} \quad \label{eq:arc_model_z}}
	{\forall u \in U}
    \end{maxi!}%
    \caption{The Arc-Based Model.}
    \label{fig:arc_model}
\end{figure}

The binary decision variable $x_{ij}^f$ determines whether vehicle $f
\in F$ traverses edge $(i,j)$. For simplicity, let $x^f(\delta^+(i)) =
\sum_{(i,j) \in \delta^+(i)} x_{i,j}^f$ and $x^f(\delta^-(i)) =
\sum_{(j,i) \in \delta^-(i)} x_{j,i}^f$ denote the sum of outgoing and
incoming flows to node $i$, respectively, where $\delta^+(i)$ and
$\delta^-(i)$ represent the sets of outgoing and incoming edges of
node $i \in N$.

The objective function, defined in equation~\eqref{eq:arc_model_obj}, 
aims at maximizing the number of trip requests served. 
\revision{Constraint~\eqref{eq:request} enforces the complete service requirement for each 
rider $u \in U$. The binary variable $z_u$ indicates whether rider $u$ 
is served ($z_u=1$) or not ($z_u=0$). The constraint links this variable 
to the service status of all pickup nodes $i$ belonging to the 
set $P_u$ (the set of pickup nodes associated with rider $u$'s requests 
$R_u$). By requiring $\sum_{f\in F}x^{f}(\delta^{+}(i))$ to equal the 
same value $z_u$ for all $i \in P_u$, this constraint ensures 
that if any one of rider $u$'s requests is served, all of their other 
requests must also be served by some vehicle in the fleet $F$.}
Constraint~\eqref{eq:same_vehicle} enforces that pickup and
drop-off services for each trip request are completed by the same
vehicle. Constraint~\eqref{eq:vehicle_flow_balance} ensures flow
balance, and constraint~\eqref{eq:arrival_time} updates the variables
$T_i^f$, arrival time of vehicle $f$ at node
$i$. Constraint~\eqref{eq:time_window} enforces time window
requirements at node $i$. Constraint~\eqref{eq:shift} specifies the
maximum working hours of each driver.
Constraint~\eqref{eq:onbard_passengers} updates the variables $Q_j^f$,
the number of riders in vehicle $f$ after visiting node
$j$. Constraint~\eqref{eq:capacity_bound_ab} specifies the bounds of
$Q_i^f$. Lastly, constraints \eqref{eq:arc_model_x} and
\eqref{eq:arc_model_z} define the domains of the decision variables.

Solutions from the arc-based model~\eqref{formulation:arc_model} provide
insights for both rider trip planning and crew shift scheduling. On
the rider side, the solutions indicate whether a rider can be served,
and if so, by which vehicle and at what time. On the driver side,
the solutions guide which riders to serve and when. In addition,
$T_0^f$ and $T_{2n+1}^f$ represent the departure and arrival times at
the depot, which can be used to determine crew shifts. While the model
formulation in Figure \ref{fig:arc_model} is nonlinear, it can be
easily linearized since the variables $x_{ij}^f$ are binary.

\subsection{A Path-Based Model}

\revisiontwo{Off-the-shell solvers such as Gurobi or CPLEX are not capable of directly solving large-scale
instances encountered in practice when using the arc-based formulation
presented in Section~\ref{sec:arc_model}. Instead, \algo{} adopts a path-based model, presented in Figure~\ref{fig:path_model}, which corresponds to a Dantzig-Wolfe decomposition of the arc-based formulation. This reformulation is well-established in the vehicle routing literature and has demonstrated strong computational performance for the classical VRP and its variants \citep{cordeau2000vrp, toth2014vehicle}.} On many
practical applications, the path-based formulation can leverage column
generation to find high-quality solutions more efficiently, as
demonstrated in prior studies (e.g.,
\cite{riley2019column,lu2022rich}).

\begin{figure}[!ht]
    \begin{maxi!}
%
	{}
%
	{\sum_{r \in R} y_r \label{eq:path_model_obj}}
%
	{\label{formulation:path_model}}
%
	{}
%
	\addConstraint
	{y_r}
	{\leq \sum_{\theta \in \Omega} \alpha_{r\theta}\lambda_\theta \quad \label{eq:path_trip_request}}
	{\forall r \in R}
	\addConstraint
	{y_r}
        {= y_{r'} \quad \label{eq:path_passenger_service}}
	{\forall u \in U, \forall r, r' \in R_u }
	\addConstraint
	{\sum_{\theta \in \Omega} \lambda_\theta}
	{\leq |F| \quad \label{eq:path_total_vehicle}}
	{}
	\addConstraint
	{\lambda_\theta}
	{\in \{0,1\} \quad \label{eq:path_route_range}}
	{\forall \theta \in \Omega}
	\addConstraint
	{y_r}
	{\in \{0,1\} \quad \label{eq:path_y_range}}
	{\forall r \in R}
    \end{maxi!}%
    \caption{The Path-Based Model.}
    \label{fig:path_model}
\end{figure}

The objective function, defined in equation~\eqref{eq:path_model_obj},
aims at maximizing the number of fulfilled trip requests. The binary
variable $y_r$ indicates whether a trip request $r$ is served. The right-hand side of constraint~\eqref{eq:path_trip_request} denotes
whether request $r$ is served by the selected routes from the set
$\Omega$. Here, $\alpha_{r\theta}$ represents whether request $r$ is
served by route $\theta \in \Omega$, and the binary variable
$\lambda_\theta$ indicates whether route $\theta$ is selected. Constraint~\eqref{eq:path_passenger_service} ensures
complete servicing of each rider requests without partial
fulfillment. Constraint~\eqref{eq:path_total_vehicle} controls the
maximum number of vehicles that can be deployed. Lastly,
constraints~\eqref{eq:path_route_range} and \eqref{eq:path_y_range}
define the domains for the decision variables $\lambda_\theta$ and
$y_r$, respectively.

Each $\theta \in \Omega$ is a feasible vehicle route that satisfies
constraints \eqref{eq:same_vehicle} to
\eqref{eq:arc_model_z}. Therefore, each vehicle route $\theta$ not
only contains trip service schedules but also provides driver shift
information. The major challenge of solving
Model~\eqref{formulation:path_model} is to find feasible routes and
construct the route set $\Omega$. The size of set $\Omega$ increases
exponentially with the number of trip requests. Hence, it is
practically impossible to enumerate all routes in the set.
Section~\ref{sec:column_generation} demonstrates how column generation
is used to iteratively add promising routes to the set $\Omega$,
avoiding the need to enumerate all possible routes.

\section{The Column Generation Algorithm}
\label{sec:column_generation}

This section introduces the column generation algorithm for solving the
\pb{}. Starting with a customized decomposition strategy, we illustrate 
three major components of the column generation in this study, including 
the master problem, the pricing subproblem, and how to find feasible 
integer solutions. \revision{Readers are referred to Appendix~\ref{app:cg_details} 
for more details.}

\subsection{Problem Decomposition Based on Driver Shifts}

As previously discussed, this goal of the \pb{} is not only to design
optimal trip service schedules but also to determine the driver shifts
that best accommodate the time-varying travel demands from the
riders. In practical operations, driver shifts typically commence at
specific times (hourly or half-hourly) for management
convenience. This observation underpins the initial step of column generation:
{\em the generation of the candidate set for driver shifts}.  It is
important to emphasize that, in the solution, multiple drivers can use
the same shifts and some shifts may not be used at all. It is the role
of column generation to determine the best driver shifts to serve as many
requests as possible.

\begin{figure}[!ht]
    \centering
    \includegraphics[width=0.9\textwidth]{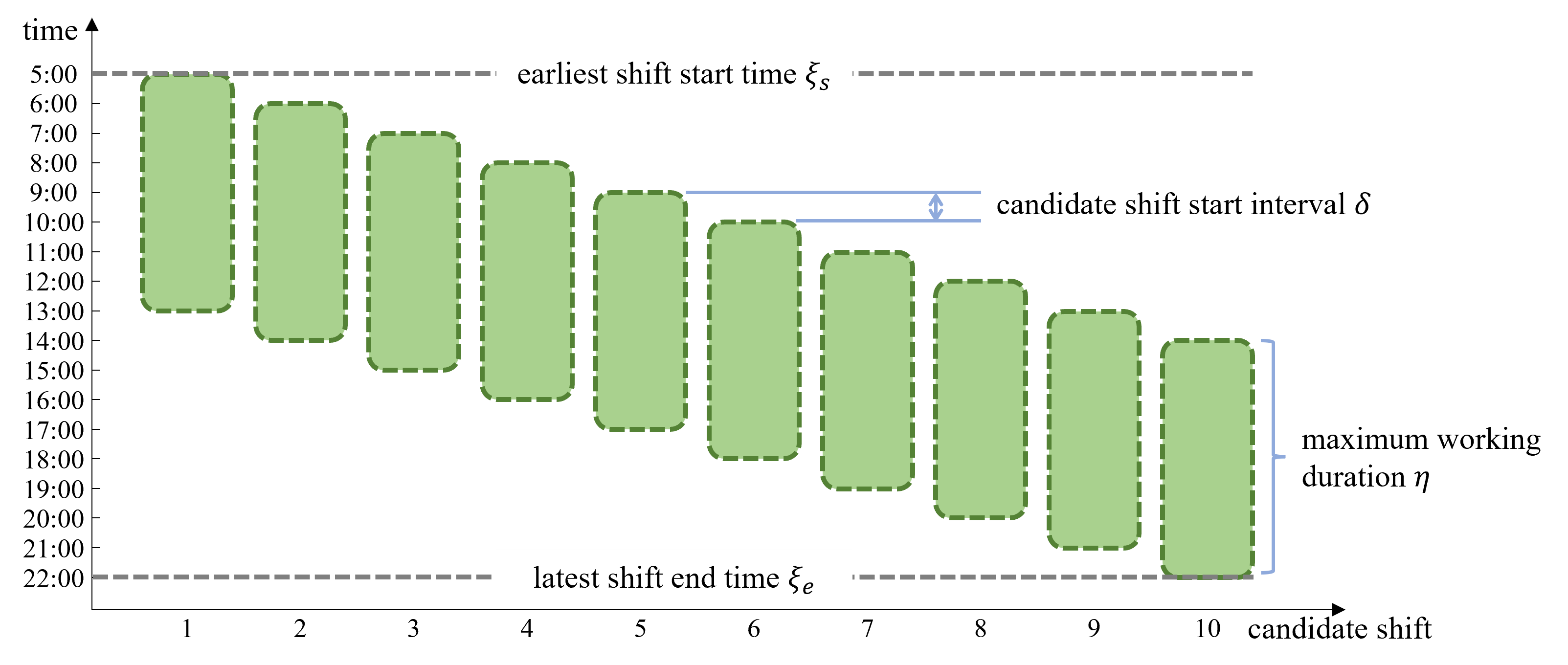}
    \caption{\revision{A Potential Candidate Set for Driver Shifts on a Particular Day.}}
    \label{fig:driver_shift_candidates}
\end{figure}

For concreteness, consider a scenario with the earliest shift start time
$\xi_s$ and latest shift end time $\xi_e$, and a maximum
driver working duration of $\eta$ hours. Assuming a time interval
$\delta$ between adjacent shift candidates, the candidate set $\Phi$
for driver shifts can be defined as
\begin{equation}
    \Phi = \{ \phi = (dr_s,dr_e) \mid dr_s = \xi_s + k\delta, dr_e = dr_s + \eta, 0 \leq k \leq \lceil \frac{(\xi_e - \xi_s - \eta)}{\delta} \rceil, k \in \mathbb{Z} \},
\end{equation}
where each element $\phi = (dr_s,dr_e) \in \Phi$ represents a driver
shift, with $dr_s$ and $dr_e$ indicating its start and end times,
respectively. Figure~\ref{fig:driver_shift_candidates} illustrates
such a candidate set for $\xi_s = 5$, $\xi_e = 22$, $\eta = 8$, and $\delta =
1$.

The Restricted Linear Master Problem (RLMP) of the column generation
algorithm is presented in Figure~\ref{fig:cg_mp}. The model uses a
given candidate set $\Phi$ of driver shifts and relaxes the
integrality of the decision variables $\lambda_\theta$ and
$y_r$. Additionally, a restricted route set $\Omega_\phi'$ is
introduced for each driver shift candidate $\phi \in \Phi$: it
contains a subset of feasible vehicle routes corresponding to driver
shift $\phi$.

\begin{figure}[!ht]
    \begin{maxi!}
%
	{}
%
	{\sum_{r \in R} y_r \label{eq:cg_mp_obj}}
%
	{\label{formulation:cg_mp}}
%
	{}
%
	\addConstraint
	{y_r}
	{\leq \sum_{\phi \in \Phi} \sum_{\theta \in \Omega_\phi'} \alpha_{r\theta}\lambda_\theta \quad \label{eq:cg_mp_trip_request}}
	{\forall r \in R}
	\addConstraint
	{y_r}
        {= y_{r'} \quad \label{eq:cg_mp_passenger_service}}
	{\forall u \in U, \forall r, r' \in R_u}
	\addConstraint
	{\sum_{\phi \in \Phi} \sum_{\theta \in \Omega_\phi'} \lambda_\theta}
	{\leq |F| \quad \label{eq:cg_mp_total_vehicle}}
	{}
	\addConstraint
	{\lambda_\theta}
	{\in [0,1] \quad \label{eq:cg_mp_route_range}}
	{\forall \theta \in \bigcup_{\phi \in \Phi} \Omega_\phi'}
	\addConstraint
	{y_r}
	{\in [0,1] \quad \label{eq:cg_mp_y_range}}
	{\forall r \in R}
    \end{maxi!}%
    \caption{The Restricted Linear Master Problem.}
    \label{fig:cg_mp}
\end{figure}

For each $\Omega_\phi'$, promising feasible routes using driver shift
$\phi$ are iteratively added by solving pricing subproblems. Pricing
subproblems corresponding to different driver shifts are completely
independent.  Therefore, in the pricing stage, the original pricing
subproblem can be decomposed to $|\Phi|$ independent pricing problems
that can be solved in parallel.

It is important to note that a finer granularity in the
driver shift set $\Phi$, determined by $\delta$, does not
substantially increase the complexity of identifying promising routes
in these subproblems. This is because each subproblem is independent
and can be processed in parallel. Furthermore, a higher number of
subproblems potentially introduces more promising routes in each
column generation iteration, which may reduce the total number of
iterations required.

\subsection{The Pricing Subproblem}
\label{sec:pricing_subproblem}

Each column generation iteration adds, for each driver shift $\phi \in
\Phi$, new promising routes in $\Omega_\phi'$ by solving the pricing
subproblem shown in Figure \ref{fig:cg_sp}.  In the objective
function~\eqref{eq:cg_sp_obj}, $\pi_r$ and $\sigma$ denote the dual
values of constraints \eqref{eq:cg_mp_trip_request} and
\eqref{eq:cg_mp_total_vehicle} after solving
Problem~\eqref{formulation:cg_mp}. 
\revisiontwo{The pricing subproblem seeks a feasible vehicle route that 
minimizes the reduced cost. This problem can be modeled as a shortest 
path problem with resource constraints, a common approach in 
column generation for vehicle routing problems \citep{feng2024branch}. The 
goal is to find a path from the starting depot to the ending depot in 
the graph $g=(N, E)$ such that the path's total reduced cost is minimized, 
while adhering to the resource constraints in \eqref{eq:same_vehicle} - \eqref{eq:arc_model_z}.

In \algo{}, this pricing subproblem is solved using dynamic programming. 
The specifics of this dynamic programming approach, including the label 
definition and extension process, are detailed in Algorithm 3 in 
Appendix~\ref{app:cg_details}.} The dynamic programming implementation does 
not solve the pricing subproblem to optimality for each driver shift, which 
can be very time-consuming. Instead, the dynamic programming process stops 
as soon as a fixed number of routes with negative reduced costs has been 
generated. Note that the early stop strategy does not affect the optimality 
of the algorithm.

\begin{figure}[!ht]
    \begin{mini!}
%
	{}
%
	{-\sum_{r \in R} \pi_r y_r -\sigma \label{eq:cg_sp_obj}}
%
	{\label{formulation:cg_sp}}
%
	{}
%
	\addConstraint
	{\eqref{eq:same_vehicle} - \eqref{eq:arc_model_z} \notag}
	{}
	{}
    \end{mini!}%
    \caption{The Pricing Subproblem For Driver Shift $\phi$.}
    \label{fig:cg_sp}
\end{figure}

Preliminary experiments indicated that, for large-scale instances,
even with the early stop strategy, it remains challenging to generate
a sufficient number of promising routes with negative costs within an
acceptable timeframe. To address this, the column generation
implements two heuristics.  First, in the dynamic programming process,
at each node $i \in N$, only a limited number of labels with the
lowest reduced costs are retained. Second, during column generation,
if the value of the objective function~\eqref{eq:cg_mp_obj} remains
unchanged for a specified number of iterations, the process is
terminated early.  It is important to note that these heuristics,
along with the one mentioned in the next subsection, may affect the
ability of the column generation algorithm in finding an optimal solution.

\subsection{Finding Feasible Integer Solutions}
\label{sec:feasible_solution}

The column generation often produces solutions that are fractional,
rendering them infeasible for the original model detailed in
formulation~\eqref{formulation:path_model}. To obtain integer
solutions, the column generation algorithm is typically embedded
within a branch-and-bound framework.  Two primary approaches are
prevalent for branching strategies: edge-based and
path-based. Edge-based branching tends to yield more balanced
subproblems, whereas path-based branching is more efficient in rapidly
identifying integer feasible solutions.

However, preliminary experiments revealed that, given the high 
complexity and large scale of the problems in this study, neither 
exact edge-based nor path-based branching strategies could produce 
feasible solutions within a reasonable time frame. Consequently, 
the column generation developed in this study employs a path-based
branching strategy complemented by a straightforward
heuristic. Specifically, if the optimal solution of the RMLP is
fractional, the variable associated with the column with the largest
fractional $\lambda_\theta$ value is fixed to 1 and a new phase of the column
generation algorithm is initiated. This process is repeated until an
integer feasible solution is obtained.

\section{Boosting Column Generation with Machine Learning}
\label{sec:gnn}

In the column generation method outlined in
Section~\ref{sec:column_generation}, the most computationally
intensive components are the pricing subproblems, which are NP-hard as
they involve finding the shortest paths with resource constraints.  To
overcome this computational challenge, \algo{} leverages machine
learning to speed up pricing subproblems. \revisiontwo{This section presents
\algo{} (Attention and Gated GNN-Informed Column Generation), a
framework that leverages machine learning to speed up and improve the
quality of column generation for real-life applications where trip
demands exhibit stable patterns over time. The framework extract insights from historical data to inform and guide the solution of new problem instances.}

\subsection{The Machine Learning Framework}
\label{subsec:framework}

\begin{figure}[!ht]
    \centering
    \includegraphics[width=0.8\textwidth]{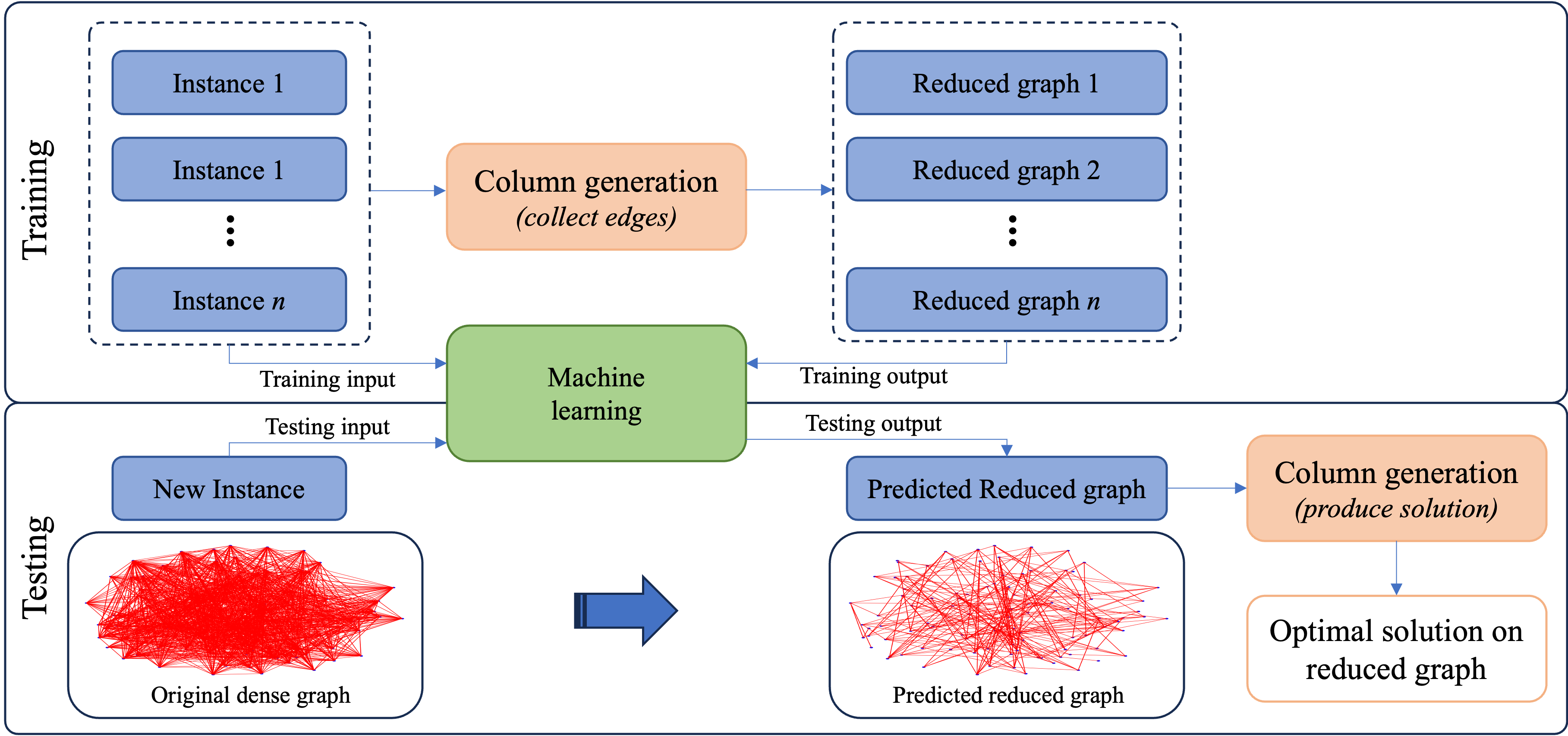}
    \caption{The Machine Learning Framework of \algo{}.}
    \label{fig:framework}
\end{figure}

The \algo{} framework is depicted in Figure~\ref{fig:framework}.  Its
key strategy is to reduce the number of edges to explore in the
pricing subproblems. In other words, before applying column
generation, \algo{} simplifies the graph of a new instance by
discarding edges unlikely to yield high-quality routes.

To determine which edges belong to high-quality or optimal routes,
\algo{} uses routes generated in the pricing subproblems of historical
instances. The training of machine learning model of \algo{} uses a dataset of
historical instances and their solutions using the the column
generation algorithm introduced in Section~\ref{sec:column_generation}. 
Edges that frequently appear across column generation iterations are 
earmarked as {\em promising}. These promising edges make it possible to 
define reduced graph for each historical instance. The historical 
reduced graphs are key input to the machine learning, as shown in the 
top part of Figure~\ref{fig:framework}.

\algo{} then trains a machine learning model in a supervised manner:
the training optimization receives as inputs the historical graphs of
each instance and the reduced graphs. In real time, \algo{} applies
the machine learning to a new instance and obtains its reduced
graph. The column generation algorithm is applied to the reduced
instance. The reduced graph, with fewer edges, allows the pricing
subproblems to be solved more quickly, greatly improving the overall
efficiency of the column generation process. In practice, the
instances come with different nodes and different edges. As a result,
\algo{} uses a Graph Neural Network (GNN) \citep{kipf2016semi} for the
learning task. For each instance, the inputs for training the GNN model
are the graph structure and encodings of the nodes and edges (see
Section \ref{section:encoding}). The GNN architecture predicts which edges
are promising for any instance represented in this way. Additionally, Section
\ref{section:promising} explores how to classify an edge as
promising in the training data.

\subsection{Overview of The GNN-Based Architecture}
\label{section:overview}

\begin{figure}[!ht]
    \centering
    \includegraphics[width=0.9\textwidth]{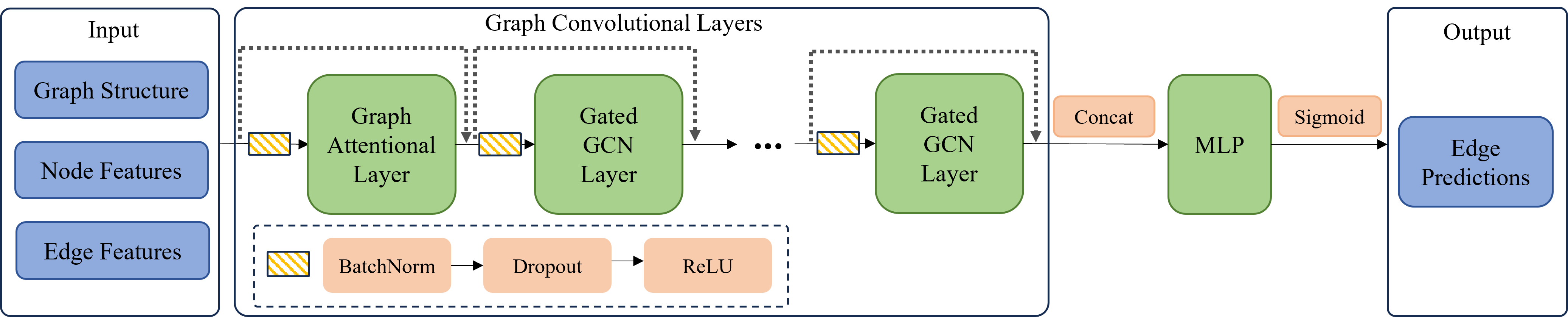}
    \caption{The GNN Architecture of \algo{}.}
    \label{fig:gnn-pipeline}
\end{figure}

The GNN architecture of \algo{} is illustrated in Figure
\ref{fig:gnn-pipeline}. 
Table
\ref{table:glossary} provides a complete glossary of terms used in the GNN model, including inputs, outputs, learnable parameters, auxiliary variables, and loss function parameters.
The GNN model is presented mathematically in Figure \ref{fig:ml_gnn}, and the training optimization problem is shown in Figure \ref{fig:ml_training}. The rest of this section introduces each concept
incrementally.

\begin{table}[!ht]
  \caption{Glossary: Inputs, Outputs, Learnable Parameters, and Loss Function Parameters.}
  \rule{\linewidth}{0.1mm}
  \begin{tabular}{lll}
  {\bf Inputs:} &                         &  \\
  & $g = (N^g, E^g)$ & the graph \\
  & $\mathbf{h}_i^{g}$         & the initial embedding of node $i$ in $g$ \\
  & $\mathbf{e}_{ij}^{g}$      & the initial embedding of edge $(i,j)$ in $g$ \\
  & $p_{ij}^g$                & the label for edge $(i,j)$ of $g$ \\
  \\
  {\bf Outputs:} &                        & \\
  & $\hat{p}_{ij}^g$                 & the prediction for edge $(i,j)$ of $g$ \\
  \\
  {\bf Learnable Parameters:} &                        & \\
  & $\mathbf{W}_{1,h}, \; \mathbf{W}_{2,h}, \; \mathbf{W}_{3,h}, \; \mathbf{a}_h$                & edge attention (head $h$) \\
  & $\mathbf{W}_{4,h}, \; \mathbf{W}_{5,h}$                             & node attention (head $h$) \\
  & $\mathbf{\Theta_{1, \ell}}, \; \mathbf{\Theta_{2, \ell}}, \; \mathbf{\Theta_{3, \ell}}, \; \mathbf{\Theta_{4, \ell}}$ & residual gated graph (layer $\ell$) \\
  & {\bf MLP} & decoder as a multi-layer perceptron \\
  & $\mathbf{\Pi}$ & the collection of all learnable parameters \\
  \\
    {\bf Auxiliary Values:} &                        & \\
    & $\mathbf{h}_i^{g,\ell}$         & transformed embedding for node $i$ in layer $\ell$ in $g$ \\
    & $\gamma_{ij}^{g,h}$  & attention of edge $(i,j)$ for head $h$ in $g$ \\
    & $\kappa_{ij}^{g,\ell}$  & gate of edge $(i,j)$ for layer $\ell$  in $g$ \\    
    
  \\
    {\bf Loss Function Parameters:} &                        & \\
    & $w^1 \; w^0$ & label class weights \\
    & $\alpha_{L}$     & L1-Norm penalty \\
  \end{tabular}
    \rule{\linewidth}{0.1mm}
    \label{table:glossary}
\end{table}

\begin{figure}[!ht]
  \rule{\linewidth}{0.1mm}
  ${\cal M}_{\mathbf{\Pi}}\left((N^g,E^g),\{\mathbf{h}_i^{g}\}_{i \in N^g},\{\mathbf{e}_{ij}^{g}\}_{(i,j) \in E^g}\right):$  \\
\begin{minipage}{\linewidth} 
  \begin{adjustwidth}{3em}{}
    \begin{flalign*}
&    \gamma_{ij}^{g,h} = 
      \frac{\exp{\left(\mathbf{a}_h^\top \, \text{LeakyReLU}(\mathbf{W}_{1,h} \, \mathbf{h}_i^g + \mathbf{W}_{2,h} \, \mathbf{h}_j^g + \mathbf{W}_{3,h} \, \mathbf{e}_{ij}^g)\right)}}{
      \sum_{k \in N(i) \cup \{i\}}\exp{\left(\mathbf{a}_h^\top \, \text{LeakyReLU}(\mathbf{W}_{1,h} \, \mathbf{h}_i^g + \mathbf{W}_{2,h} \, \mathbf{h}_k^g  + \mathbf{W}_{3,h} \, \mathbf{e}_{ik}^g)\right)}} & (i,j) \in E^g, \, h \in H \\
      &    \mathbf{h}_i^{g,1} = \mathop{\Vert}\limits_{h=1}^{H} \gamma_{ii}^{g,h} \mathbf{W}_{4,h} \mathbf{h}_i^{g}  + 
\sum_{j \in N(i)} \gamma_{ij}^{g,h} \mathbf{W}_{5,h} \mathbf{h}_j^g & i \in N^g \\
      & \kappa_{ij}^{g,\ell} = \sigma(\mathbf{\Theta_{3, \ell}} \,  \mathbf{h}_i^{g,\ell -1} + \mathbf{\Theta_{4, \ell}} \,  \mathbf{h}_j^{g,\ell -1}) & (i,j) \in E^g, \ell \in 2 \ldots L \\
      & \mathbf{h}_i^{g,\ell} = \mathbf{\Theta_{1, \ell}} \, \mathbf{h}_i^{g,\ell -1} + \sum_{j \in N(i)} \kappa_{ij}^{g,\ell} \odot \mathbf{\Theta_{2, \ell}} \, \mathbf{h}_j^{g,\ell -1}  & i \in N^g, \ell \in 2 \ldots L \\
& \hat{p}_{ij}^g = \sigma \left(\text{{\bf MLP}}(\{\mathbf{h}_i^{g,L} \mathbin\Vert \mathbf{h}_j^{g,L}\}) \right) & (i,j) \in E^g \\
      &  \mbox{{\bf return}} \{\hat{p}_{ij}^g\}_{(i, j) \in E^g} \\
    \end{flalign*}
 \end{adjustwidth}
\end{minipage}
  \rule{\linewidth}{0.1mm}
    \caption{The Machine Learning GNN Model for Predicting Promising Edges.}
    \label{fig:ml_gnn}
\end{figure}

\begin{figure}[!ht]
  \rule{\linewidth}{0.1mm}  
    \begin{mini!}
%
	{\mathbf{\Pi}}
%
	{\mathlarger{\sum_{m=1}^M} \ \sum_{(i,j) \in E^{g_m}} \left(w^1 \, p_{ij}^{g_m} \, \log(\hat{p}_{ij}^{g_m}) + w^0 \, (1 - p_{ij}^{g_m}) \, \log(1 - \hat{p}_{ij}^{g_m}) \right) + \alpha_{L} \ \sum_{(i,j) \in E^{g_m}} \, \hat{p}_{ij}^{g_m} \notag}
%
	{\label{formulation:gnn_training}}
%
	{}
%
	\addConstraint
	{\{\hat{p}^{g_m}_{ij}\}_{(i,j) \in E^{g_m}} \, = \, }
	{{\cal M}_{\mathbf{\Pi}}\left((N^{g_m},E^{g_m}),\{\mathbf{h}_i^{g_m}\}_{i \in N^{g_m}},\{\mathbf{e}_{ij}^{g_m}\}_{(i,j) \in E^{g_m}}\right)}
	  {\quad \forall m \in M \notag}
    \end{mini!}
      \rule{\linewidth}{0.1mm}
    \caption{The Optimization Model for Training the GNN Model.}
    \label{fig:ml_training}
\end{figure}

As mentioned earlier, the GNN receives, as
input, a graph structure, the node encodings, and the edge encodings.
The GNN then applies several graph convolutional layers, starting with a
graph attention layer before a series of gated graph convolutional layers. These
layers transform the node encodings. The graph convolutional layers are
followed by a multi-layer perceptron that takes, as inputs, the
transformed encodings of every edge in the graph and produces an
output that is then run through a sigmoid function to estimate the
likehood of an edge to be promising. More formally, the GNN
architecture defines a parametric function
\[
  {\cal M}_{\mathbf{\Pi}}\left((N^g,E^g),\{\mathbf{h}^g_i\}_{i \in N^g},\{\mathbf{e}^g_{ij}\}_{(i,j) \in E^g}\right)
  \]
which receives, as inputs, a graph structure $(N^g,E^g)$, a
node encoding $\{\mathbf{h}^g_i\}_{i \in N^g}$, and an edge
encoding $\{\mathbf{e}^g_{ij}\}_{(i,j) \in E^g}$. It returns the
likelihood $\hat{p}^g_{ij}$ that an edge be promising for the column
generation for each edge $(i,j) \in E^g$. The function is
parametrized by a set of learnable parameters $\mathbf{\Pi}$.

The GNN architecture uses the graph structure as an input, so that it
can be applied to graphs of different sizes and shapes. Moreover, node
and edge encodings make it possible to apply the architecture for
instances with nodes and edges that have not been seen before. The GNN
architecture transforms a node encoding through multiple layers, using
its encoding and the encodings of its neighbors. The learnable
parameters do not depend on the specific graph structure of each
instance: they are {\em shared} by all instances and only apply to the
node and edge embeddings. GNNs  are especially
attractive in the setting of this paper, since \pb{} instances may
have different sets of nodes and different sets of edges. Indeed, the
number of requests can fluctuate significantly from day to day,
resulting in graphs of different sizes.

\subsection{Feature Encoding}
\label{section:encoding}

For each graph instance $(N^g,E^g)$, \algo{} extracts
features to characterize each node $i \in N^g$. Specifically,
the feature vector $\mathbf{h}^g_i$ for node $i$ include its geographic
coordinates (latitude and longitude), its time window, and a
categorical indicator showing whether the node is a pickup, drop-off,
or depot.  In addition, \algo{} considers features for each edge $(i,
j) \in E^g$. The edge feature vector $\mathbf{e}^g_{ij}$
includes the travel time from node $i$ to node $j$ and a binary
indicator of whether $(i ,j)$ directly connects the origin and
destination of the same trip. To ensure consistency across different
input instances, \algo{} applies a min-max scaler to normalize each
feature to the range $[0, 1]$. \revisiontwo{This feature selection aligns with standard practices in the literature, following similar approaches in machine learning applications for routing problems \citep{joshi2019efficient, yuan2022neural,  morabit2023machine}. Appendix \ref{sec:feature_importance} provides a detailed analysis of the relative importance of the proposed features.}

\subsection{The GNN Architecture}
\label{section:gnn-layers}

The section details the GNN architecture for an input
$
\left((N^g,E^g),\{\mathbf{h}_i^{g}\}_{i \in N^g},\{\mathbf{e}_{ij}^{g}\}_{(i,j) \in E^g}\right).
$

\subsubsection{The Graph Convolutional Layers}
\label{section:convolution}

The node embeddings are first updated through one layer of the
multi-head graph attenional operator (GATConv), the core of the Graph
Attention Network (GAT) architecture from \citep{velivckovic2017graph,
  brody2021attentive}. In contrast to classical graph convolutional
networks which employ equal-weight neighborhood aggregation, the
multi-head attention mechanism allows for the assignment of different
weights to different neighbor nodes. Taking the edge encodings into
consideration, \algo{} computes the attention weight
$\gamma_{ij}^{g,h}$ between a pair of nodes $i$ and $j$ for a head $h$
as follows:
\begin{equation}
\gamma_{ij}^{g,h} = 
      \frac{\exp{\left(\mathbf{a}_h^\top \, \text{LeakyReLU}(\mathbf{W}_{1,h} \, \mathbf{h}_i^g + \mathbf{W}_{2,h} \, \mathbf{h}_j^g + \mathbf{W}_{3,h} \, \mathbf{e}_{ij}^g)\right)}}{
      \sum_{k \in N(i) \cup \{i\}}\exp{\left(\mathbf{a}_h^\top \, \text{LeakyReLU}(\mathbf{W}_{1,h} \, \mathbf{h}_i^g + \mathbf{W}_{2,h} \, \mathbf{h}_k^g  + \mathbf{W}_{3,h} \, \mathbf{e}_{ik}^g)\right)}} 
\end{equation}
\begin{equation}
    \text{LeakyReLU}(x) = 
    \begin{cases} 
    x & \text{if } x > 0, \\
    0.2 x & \text{if } x \leq 0.
    \end{cases},
\end{equation}
\noindent
where $N(i)$ corresponds to the neighbors of node $i$, while 
$\mathbf{a}$, $\mathbf{W}_1$, $\mathbf{W}_2$, and $\mathbf{W}_3$ are learneable
parameters. Intuitively, $\gamma_{ij}$ measures the significance of node $j$ for node $i$.

\begin{figure}[!ht]
    \centering
    \begin{minipage}{0.35\textwidth}
        \centering
        \includegraphics[width=\linewidth]{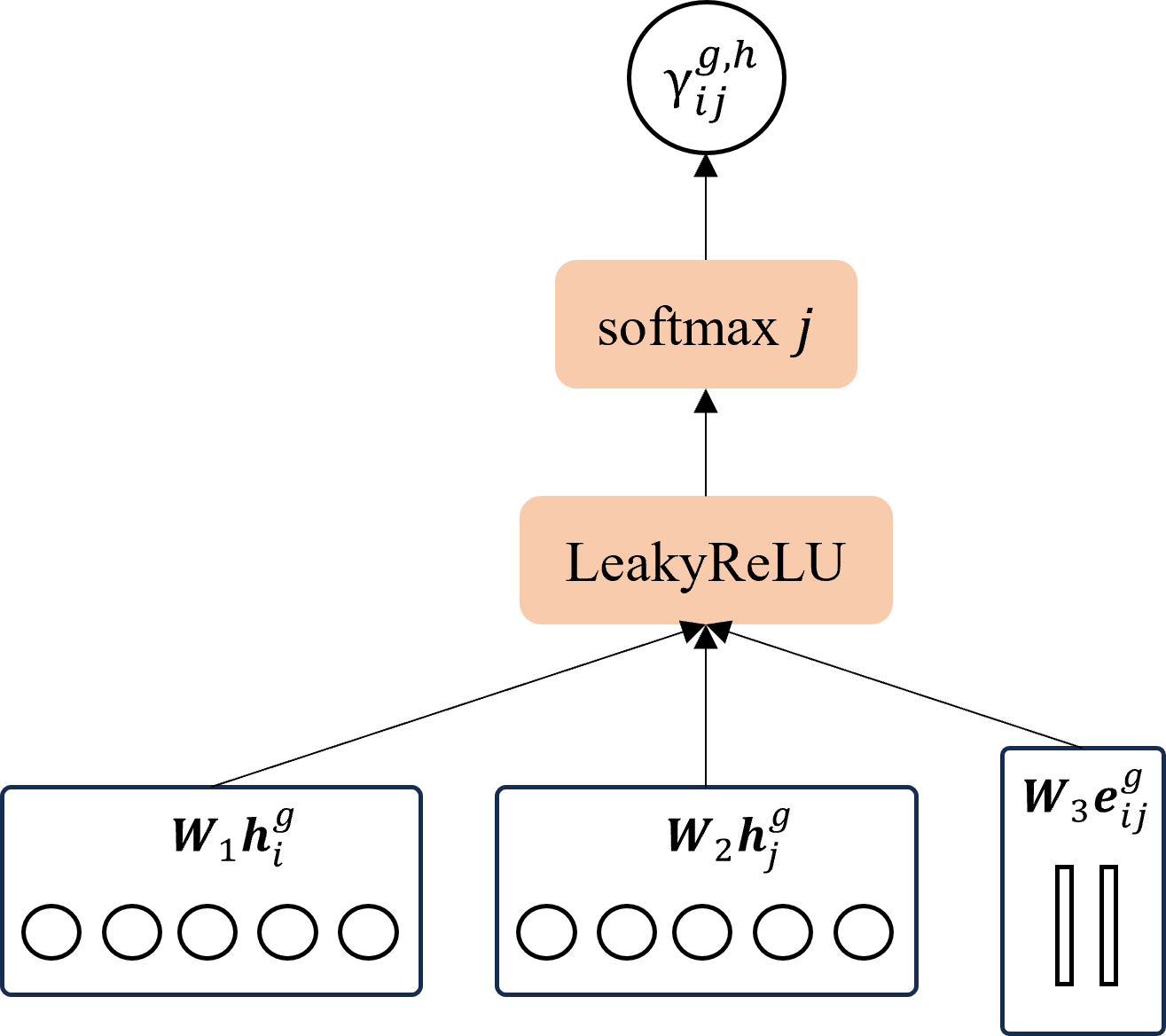}
        \caption{Visualization of the Attention Mechanism.}
        \label{fig:attention-coefficent}
    \end{minipage}\hfill
    \begin{minipage}{0.5\textwidth}
        \centering
        \includegraphics[width=\linewidth]{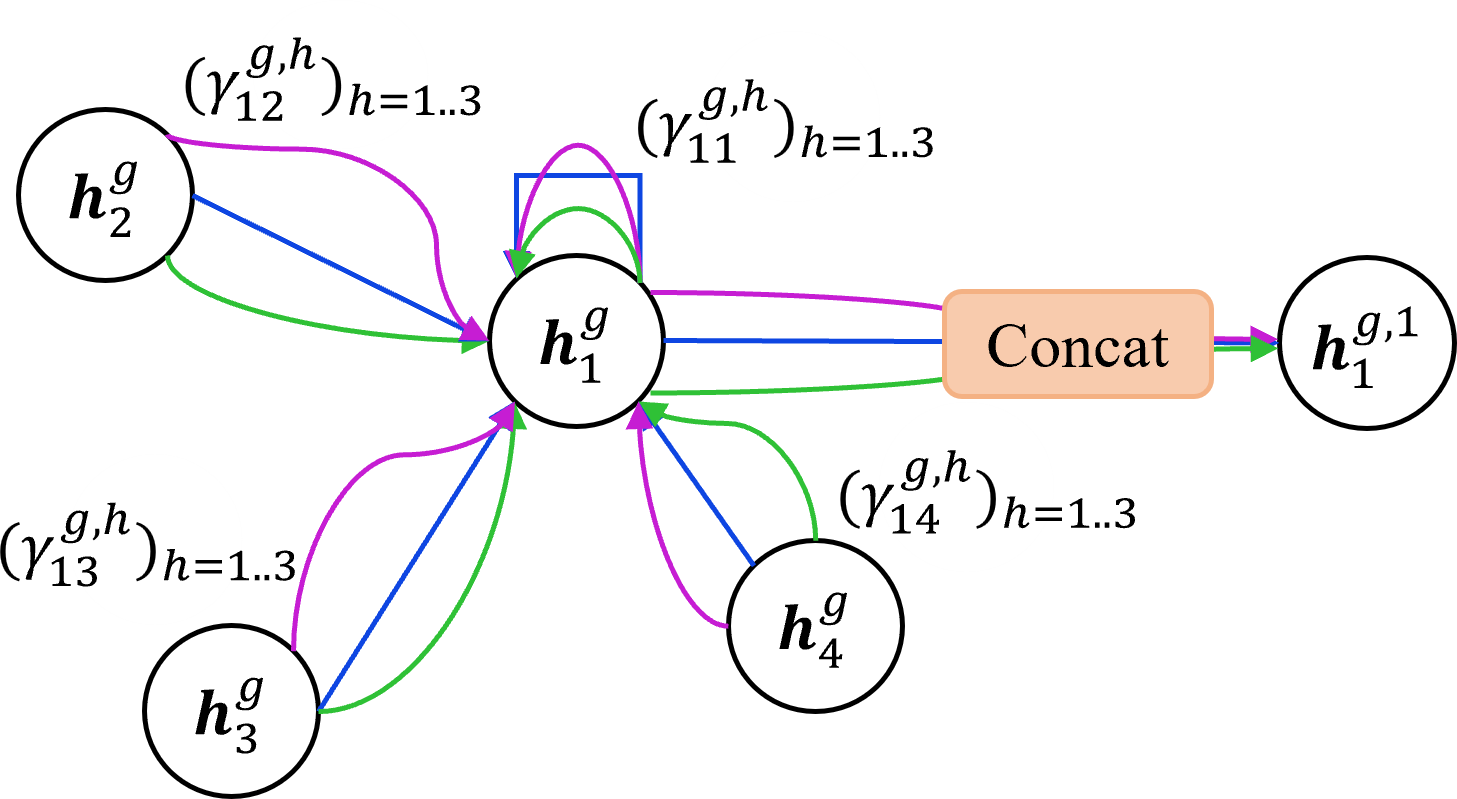}
        \caption{Visualization of Multi-Head Attention (with $H$ = 3 heads) by Node 1 on its Neighborhood; different arrow colors represent different attention heads; adapted from \cite{velivckovic2017graph}. }
        \label{fig:multi-head}
    \end{minipage} \hfill
     \begin{minipage}{0.5\textwidth}
        \centering
        \includegraphics[width=\linewidth]{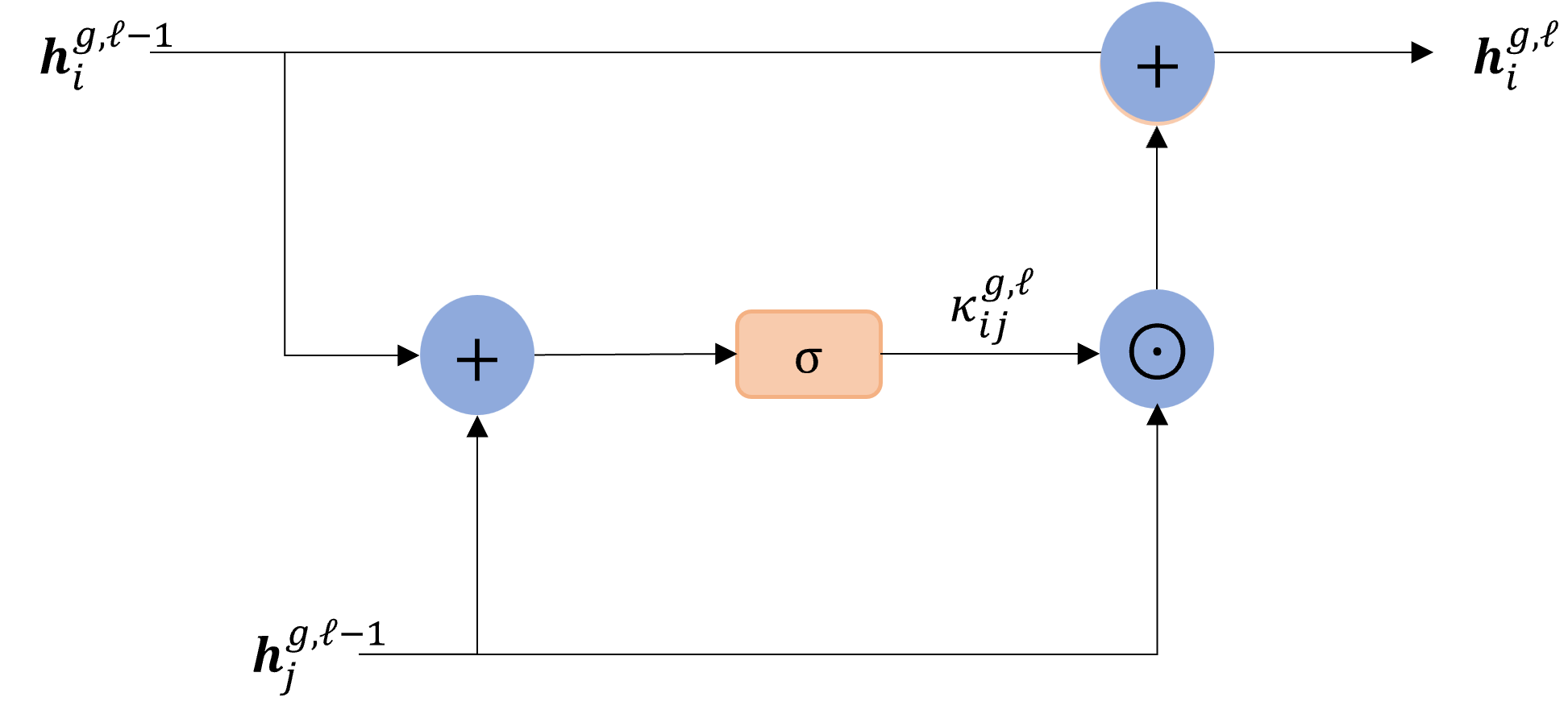}
        \caption{Visualization of the Gating Mechanism.}
        \label{fig:gate}
    \end{minipage}
\end{figure}

After obtaining the attention weights for each edge, \algo{} updates
the node embeddings of each node $i$ by calculating the weighted sum
of the transformed features of its neighbors and itself, i.e.,
\begin{equation} \label{eq: attention}
\gamma_{ii}^{g,h} \, \mathbf{W}_{4,h} \, \mathbf{h}_i^g  + \sum_{j \in N(i)} \gamma_{ij}^{g,h} \, \mathbf{W}_{5,h} \, \mathbf{h}_j^g 
\end{equation}
Since the implementation uses a multi-head attention with $H > 1$
attention heads, the node embeddings are updated by concatenating $H$
independent attention mechanisms following Equation (\ref{eq:
  attention}):
\begin{equation}
\mathbf{h}_i^{g,1} = \mathop{\Vert}\limits_{h=1}^{H} \gamma_{ii}^{g,h} \mathbf{W}_{4,h} \mathbf{h}_i^{g}  + 
\sum_{j \in N(i)} \gamma_{ij}^{g,h} \mathbf{W}_{5,h} \mathbf{h}_j^g 
\end{equation}
where $\mathop{\Vert}$ denotes the concatenation operator, while
$\gamma_{ij}^{g,h}$ represents the attention weights for the
$h$-th attention mechanism and edge $(i, j)$. The matrices
$\mathbf{W}_{4,h}$ and $\mathbf{W}_{5,h}$ are the corresponding
learnable parameters. Figures \ref{fig:attention-coefficent} and
\ref{fig:multi-head} provide visualizations of the process occurring
within the graph attentional layer.

Subsequently, the node embeddings are fed into multiple layers of the
Residual Gated Graph ConvNets (GCN), which has been used in several
related studies \citep{joshi2019efficient, yuan2022neural}. The number
of Residual Gated GCN layers is considered as a hyperparameter, which
is fine-tuned during the experiments. For each layer $2 \leq \ell \leq
L$, the node embeddings are updated following a gating mechanism
described in \citep{bresson2017residual}:
\begin{equation}
\mathbf{h}_i^{g,\ell} = \mathbf{\Theta}_{1, \ell} \, \mathbf{h}_i^{g,\ell -1} + \sum_{j \in N(i)} \kappa_{ij}^{g,\ell} \odot \mathbf{\Theta}_{2, \ell} \, \mathbf{h}_j^{g,\ell -1}
\end{equation}
with the gate $\kappa_{ij}^{g,\ell}$ defined as:
\begin{equation}
\kappa_{ij}^{g,\ell} = \sigma(\mathbf{\Theta}_{3, \ell} \,  \mathbf{h}_i^{g,\ell -1} + \mathbf{\Theta}_{4, \ell} \,  \mathbf{h}_j^{g,\ell -1}) 
\end{equation}
where $\odot$ denotes the element-wise multiplication operator. The
matrices $\mathbf{\Theta}_{1, \ell}, \mathbf{\Theta}_{2, \ell},
\mathbf{\Theta}_{3, \ell}, \mathbf{\Theta}_{4, \ell}$ are learneable
parameters for the $\ell$-th gated GCN layer, and $\sigma(\cdot)$ represents
the Sigmoid function. The gating mechanism is illustrated in Figure
\ref{fig:gate}.

The model performance was enhanced by applying established methods to
the node embeddings before they enter each graph convolutional
layer. Batch normalization and dropouts stabilize training and prevent
overfitting. Moreover, the GNN also incorporates skip (residual)
connections to amplify the influence of the embeddings generated in
the initial layers. This approach has shown empirical success in
improving node differentiation within the model.

\subsubsection{The Decoder}

After the graph convolutional layers generate the final node
embeddings, \algo{} converts them into edge embeddings by
concatenating the embeddings of the two nodes that form each edge. The
concatenated embeddings are then passed through a multi-layer
perceptron (MLP) consisting of fully connected layers followed by ReLU
activations. Finally, a Sigmoid function is applied to the MLP outputs
to ensure they fall within $[0,1]$.  The result is an estimated
probability score $\hat{p}_{ij}^g$ for each edge, representing the
likelihood that edge $(i,j)$ belongs to the set of promising
edges. The process is captured via the following equation:
\begin{equation}
\hat{p}_{ij}^g = \sigma \left(\text{{\bf MLP}}(\{\mathbf{h}_i^{g,L} \mathbin\Vert \mathbf{h}_j^{g,L}\}) \right).
\end{equation}

\subsection{The Training Loss Function}
\label{section:loss}

The training optimization is depicted in Figure
\ref{fig:ml_training}. It receives, as input, a collection of $M$
instances with the associated labels for each edge. The training
optimization uses the GNN presented earlier to predict the likelihood
of an edge to be promising. The training uses binary cross-entropy as
its loss function during training, which is a standard choice for most
binary classification tasks, and minimizes empiral risk, i.e.,
\[
\mathlarger{\sum_{m=1}^M} \ \sum_{(i,j) \in E^{g_m}} \left(w^1 \, p_{ij}^{g_m} \, \log(\hat{p}_{ij}^{g_m}) + w^0 \, (1 - p_{ij}^{g_m}) \, \log(1 - \hat{p}_{ij}^{g_m}) \right), 
\]
where instance $m$ is specified by $((N^{g_m},E^{g_m}),\{\mathbf{h}_i^{g_m}\}_{i \in N^{g_m}},\{\mathbf{e}_{ij}^{g_m}\}_{(i,j) \in E^{g_m}})$ and the labels $\{p_{ij}^{g_m}\}_{(i,j) \in E^{g_m}}$. 
To further enhance the model performance, \algo{} incorporates an
L1 regularization term in the loss function
\[
\alpha_{L} \ \sum_{(i,j) \in E^{g_m}} \, \hat{p}_{ij}^{g_m},
\]
encouraging sparsity in the model's output probabilities.

\subsection{Labeling Historical Data}
\label{section:promising}

It remains to specify how to label historical data for training the
machine learning architecture. Indeed, the effectiveness of machine
learning model training hinges on the criteria for selecting edges to
construct the reduced graph for each training instance. In general, it
is not sufficient to select edges from the paths chosen in the optimal
solutions. These sets may not allow the subsequent column generation
procedure to find feasible solutions to unseen instances or find
solutions of high quality.  Moreover, the complexity of \pb{}
led to a number of heuristic decisions in the column generation
procedure that may negatively impact the learning process. For these
reasons, the promising edges used in the labeling process of \algo{}
come from five classes:
\begin{enumerate}
    \item RGTrn-A: these are the edges from all paths explored during 
        column generation iterations, i.e., paths in 
    $\bigcup_{\phi \in \Phi} \Omega_\phi'$.
    \item RGTrn-U: these are the edges selected in a solution (path) of the RLMP, i.e., they have a non-zero value in an RLMP.
    \item RGTrn-U$x$: The top $x \in \{80, 50, 30\}$ percent of edges
      in RGTrn-U, ranked by number of times they were used in all
      RLMP solutions.
\end{enumerate}
Note that the label is 1 for a promising edge (i.e., an edge in the
selected set); it is zero otherwise.

\subsection{Using Reduced Graph in Column Generation}

\revision{During inference, the GNN model described in Figure~\ref{fig:ml_gnn} generates probability estimates $\hat{p}_{ij}^{g} \in [0, 1]$ for each edge $(i,j)$ in an unseen instance graph $g=(N^g, E^g)$. These estimates represent the likelihood that an edge belongs in the reduced graph used for column generation pricing subproblems. However, to create a practical reduced graph, these continuous probabilities must be converted to binary decisions (i.e., include or exclude the edge).

To systematically control the reduced graph's size, \algo{} employs a rank-based thresholding approach:

\begin{enumerate}
    \item First, all edges in the original graph are ranked in descending order according to their predicted probability values $\hat{p}_{ij}^{g}$.
    \item Only the top $\tau$ percent of these sorted edges are retained for the reduced graph.
    \item To ensure solution feasibility, certain critical edges are always preserved regardless of their predicted probabilities. Specifically, for each trip request $r \in R$ with pickup node $i$ and drop-off node $n+i$, the edges connecting the origin depot to pickup node, pickup node to drop-off node, and drop-off node to destination depot $\{(0, i), (i, n+i), (n+i, 2n+1)\}$ are retained.
\end{enumerate}

This thresholding process can be formally represented as an operator $\mathcal{T}_{\tau}$ that that takes the set of predicted probabilities and the original edge set as inputs, returning the reduced edge set $E_{\text{reduced}-\tau}^g$:
\begin{equation}
E_{\text{reduced}-\tau}^g = \mathcal{T}_{\tau}\biggl({\cal M}_{\mathbf{\Pi}}\left((N^g,E^g),\{\mathbf{h}_i^{g}\}_{i \in N^g},\{\mathbf{e}_{ij}^{g}\}_{(i,j) \in E^g}\right), E^g \biggr).
\end{equation}

The resulting reduced graph 
$g_{\text{reduced}-\tau} = (N^g, E_{\text{reduced}-\tau}^g)$ is then 
used in the column generation process. Specifically, in Algorithm~\ref{alg:outer},
replace graph $g=(N^g,E^g)$ with reduced graph 
$g_{\text{reduced}-\tau} = (N^g, E_{\text{reduced}-\tau}^g)$ as an input.

The experiments in this paper explore the trade-off between computation time and solution quality by varying the threshold $\tau$ between 2 and
30. For each value of $\tau$, the corresponding class of reduced test graphs is referred to as
RGTst-$\tau$.
}

\section{Numerical Experiments}
\label{sec:experiment}

This section presents the experimental results. Section
\ref{subsec:dataset} describes the dataset, Section
\ref{subsec:baselines} presents the baselines, Section
\ref{section:gnn_eval} reports the evaluation of the GNN, Section
\ref{section:cg_eval} compares \algo{} with the baselines under
different configurations, and Section \ref{section:sensitivity}
describes the sensitivity analysis.

\subsection{Dataset Description}
\label{subsec:dataset}

\algo{} is evaluated using a real-world dataset derived from the
Paratransit service in Chatham County, Georgia, U.S. This service
caters to individuals with disabilities who are unable to use the
regular public transportation system, offering them a
reservation-based travel option. Riders are required to schedule their
trips one day in advance, with the system ceasing to accept requests
at 4 pm each day. At this cutoff time, the optimization of the driver
shifts and trip schedules for the subsequent day begins. Due to
constraints in the availability of drivers and vehicles, not all
requests can be accommodated in general. The dataset encompasses
information on daily trip requests spanning from January 2014 to
December 2019. Each request includes the rider's origin, destination,
and preferred pick-up and drop-off times. 

\begin{figure}[!ht]
    \centering
    \includegraphics[width=0.9\textwidth]{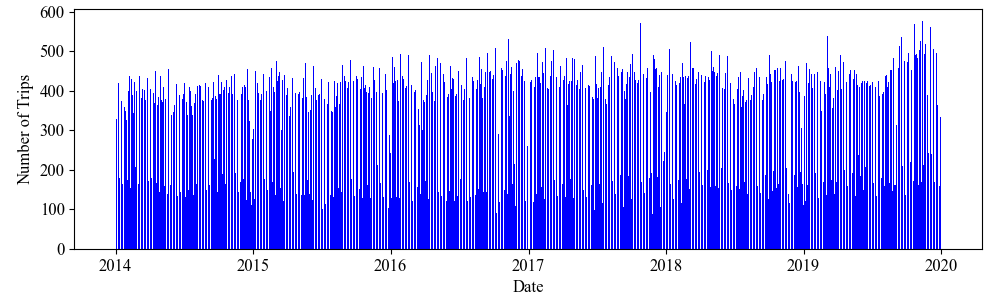}
    \caption{Daily Number of Trip Requests from January 2014 to December 2019.}
    \label{fig:number_of_trip_requests}
\end{figure}

\begin{figure}[!ht]
    \centering
    \includegraphics[width=0.9\textwidth]{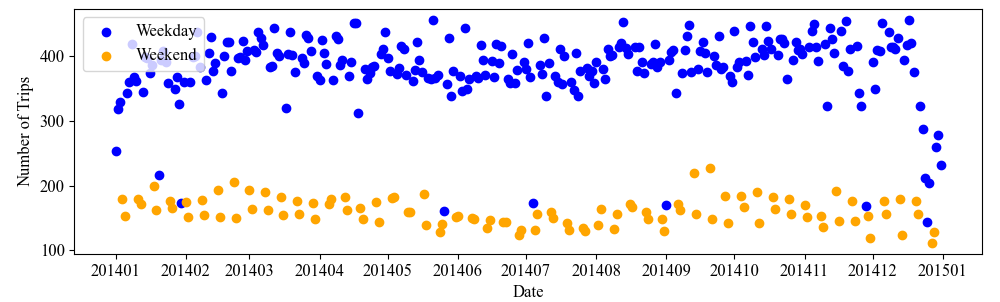}
    \caption{Daily Number of Trip Requests in 2014.}
    \label{fig:week_number_of_trip_requests}
\end{figure}

As illustrated in Figure~\ref{fig:number_of_trip_requests}, the demand
for trips is relatively stable year-over-year, an observation that is
critical for the applicability of \algo{}.
Figure~\ref{fig:week_number_of_trip_requests} indicates that weekdays
experience a higher demand for trips compared to weekends and
holidays. To ensure consistency in the demand pattern, the experiments
exclusively utilize weekday trip request data. The training of the
graph neural network described in Section~\ref{sec:gnn} is conducted
using data from workdays between January 2014 and November 2019. Data
from December 2019 serves as the basis for evaluating performance.

For the purposes of this study, the number of drivers allocated each
day is determined by the volume of trip requests for that day,
calculated as follows:
\begin{equation}
    n_d = \left\lceil \frac{n_r}{\omega} \right\rceil,
\end{equation}
where $n_d$ denotes the number of drivers, $n_r$ represents the number
of trip requests, and $\omega$ is the average number of trips a driver
can manage in a day, which is set to 16 in this study.

\subsection{Column Generation Approaches}
\label{subsec:baselines}

The following three column generation approaches are compared in the evaluation. 
\begin{enumerate}
    \item Baseline-CG is the conventional column generation approach, 
which applies column generation described in Section \ref{sec:column_generation} directly on original dense 
graphs without graph reduction. 
    \item RF-CG uses a Random Forest-based method for graph reduction, as outlined in \citet{morabit2023machine}. This approach has proven effective in a related Vehicle Routing Problem with Time Windows.
    \item \algo{} is the proposed approach which leverages a GNN-based model for graph reduction.
\end{enumerate}

By default, both RF-CG
and \algo{} use the same configuration for graph reduction training
and testing, i.e., RGTrn-U50 and RGTst-10, respectively.

\subsection{Evaluation of the Graph Neural Network}
\label{section:gnn_eval}

This section evaluates the performance of the proposed GNN model on
the real-world dataset described in Section \ref{subsec:dataset}. The
model training was conducted on a single Nvidia Tensor Core H100 GPU
on the PACE cluster \citep{PACE2017}. The implementation of the graph
convolutional layers was carried out using version 2.4.0 of the
PyTorch Geometric library in Python 3.9. The loss function is minimized via gradient
descent, using the Adam optimizer with weight decay as described in
\cite{kingma2014adam}. Each data point in the dataset corresponds to a
single graph. For graph-level mini-batching, the training follows the
method outlined in \cite{pytorch_geometric_batching}, where adjacency
matrices for the same mini-batch are stacked diagonally, and the
node-level features are concatenated along the node dimension.

Given the disparity between class distributions, with the positive
class accounting for merely 6\% of the labels, the evaluation metrics
are carefully selected to reflect a balanced view of the model's
performance. In line with \cite{morabit2023machine}, the results
utilize Recall, Specificity, and Balanced Accuracy as primary metrics.
{\em Recall}, i.e., the true positive rate, assesses the model
effectiveness in accurately identifying all instances of the positive class. {\em
  Specificity}, i.e., the true negative rate, evaluates how precisely
the model detects negatives in the majority class. As the aim is to
maximize both these metrics, {\em Balanced Accuracy}, i.e., the
average of Recall and Specificity, is also included to offer a
comprehensive metric reflecting the model's overall accuracy for both
classes.

The dataset, excluding the test period, spans from January 2014 to
November 2019 and includes a total of 1,539 instances. The dataset is
partitioned in an 80-20 split for training and validation sets. An
adaptive learning strategy is applied to optimize the training. If the
validation loss showed no improvement over several epochs, indicating
a learning plateau, the learning rate is lowered to nudge the model
towards better performance. In cases where this adjustment did not
lead to any further gains, and to prevent the model from overfitting,
early stopping is initiated. This step was crucial to retain the
model’s ability to generalize to new data. \revision{Additionally, the
hyperparameters of the GNN model are tuned to maximize balanced
accuracy on the validation set using Optuna \citep{akiba2019optuna} with 100 trials. Table \ref{tab:model_config} presents the hyperparameter search space and the best configuration identified during the tuning process.}
\begin{table}[!ht]
    \caption{\revision{GNN Hyperparameter Configuration: Search Space and Best Values Found}}
    \centering
    \begin{tabular}{@{}lcc@{}}
    \toprule
    \textbf{Hyperparameter} & \textbf{Search Space} & \textbf{Best Value Found} \\
    \midrule
    Dropout Rate & $\{0.30, 0.35, \ldots, 0.70\}$ & $0.35$ \\
    Hidden Dimensions & $\{64, 128, 256, 512\}$ & $256$ \\
    L1 Regularization ($\lambda_{L}$) & $[1 \times 10^{-8}, 1 \times 10^{-2}]$  & $4 \times 10^{-8}$  \\
    Learning Rate & $[1 \times 10^{-5}, 1 \times 10^{-2}]$  & $4 \times 10^{-3}$ \\
    Number of Attention Heads & $\{4, 8, 16\}$ & $8$ \\
    Res-Gated GCN Layers & $[2, 8]$ & $6$ \\
    Weight Decay & $[1 \times 10^{-8}, 1 \times 10^{-5}]$  & $6.8 \times 10^{-5}$ \\
    Mini-batch Size & Fixed & $4$ \\
    Class Weight Ratio ($w_{+}:w_{-}$) & Fixed & $15.65:1$ \\
    \bottomrule
    \end{tabular}
    \label{tab:model_config}
\end{table}

The performance of the GNN model is compared to a Random Forest (RF)
classification model described in \citet{morabit2023machine}. Unlike
GNN, which leverages the entire graph structure, the RF model operates
only at the edge level. RF learns the promising probability of each
edge individually based on the features of the edge and its two
endpoints. Since each training data point for RF is an edge, 30 days
from the period between January 2014 and November 2019 were sampled
to keep the training set size manageable; all edges from each of these
days are included. The features used for RF are aligned with those in
the GNN model. In addition, the balanced class weights are adopted,
and the hyperparameters are tuned similarly as in
\citet{morabit2023machine} to minimize the balanced accuracy.

The test set metrics, shown in Table \ref{tab:test_data_metrics},
highlight the strong performance of the GNN model. With a recall of
90.2\% (compared to 81.1\% for RF), the GNN accurately identifies the
majority of promising edges. Additionally, its specificity of 86.4\%
(compared to 77.8\% for RF) demonstrates that the high recall does not
come at the expense of incorrectly classifying negative instances.
The balanced accuracy of 88.3\% for the GNN is nearly 9\% higher than
for the RF model, further indicating that GNN is a more effective
classifier overall. The results demonstrate that the GNN benefits from
learning from the global graph topology, rather than relying solely on
local connectivity information.

\begin{table}[!ht]
    \caption{Accuracy Metrics on the Test Set.}
    \centering
    \begin{tabular}{lcc}
    \hline
    \textbf{Metric} & \textbf{GNN (\%)} & \textbf{RF (\%)} \\
    \hline
    Recall & 90.2 & 81.1 \\
    Specificity & 86.4 & 77.8\\
    Balanced Accuracy & 88.3 & 79.4\\
    \hline
    \end{tabular}
    \label{tab:test_data_metrics}
\end{table}

\begin{figure}[!ht]
    \centering
    \includegraphics[width=0.5\textwidth]{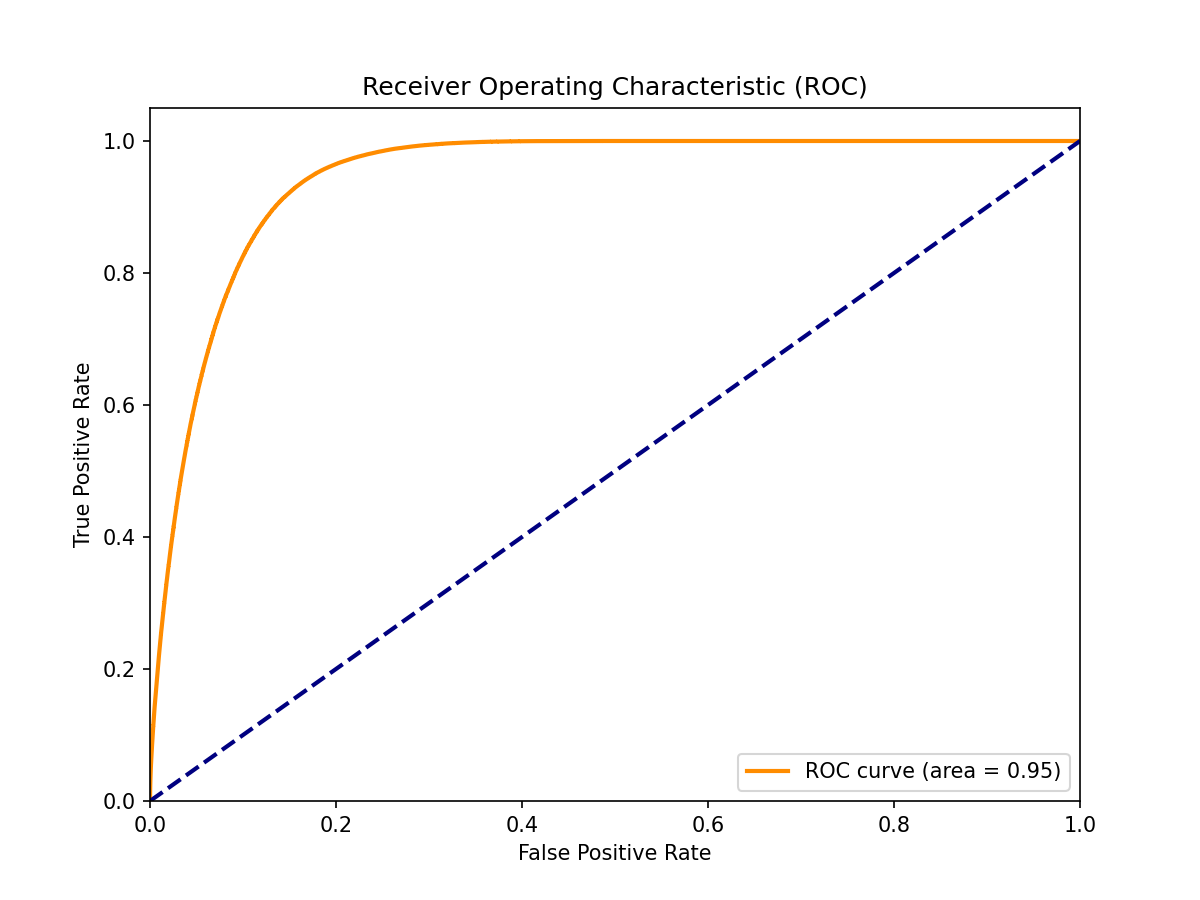}
        \caption{The ROC Curve for GNN Evaluated on the Test Set.}
        \label{fig:final_test_roc_curve}
\end{figure}

The Receiver Operating Characteristic (ROC) curve for GNN, shown in
Figure \ref{fig:final_test_roc_curve}, corroborates these
findings. With an AUC of 0.95, the curve articulates the model's
substantial discriminative power, signifying a high true positive rate
at various threshold levels while maintaining a low false positive
rate. This high AUC is particularly telling of the model's proficiency
in distinguishing between the two classes. The ROC curve's proximity
to the upper left corner are indicative of an almost ideal classifier,
striking an effective balance between sensitivity and
specificity. This is particularly notable given the challenge of
maintaining high sensitivity in this highly imbalanced dataset without
compromising specificity.

\subsection{Performance Analysis of the Column Generation}
\label{section:cg_eval}

This section compares \algo{} with Baseline-CG and RF-CG. The analysis
focuses on weekday data to account for the notable differences in trip
request patterns between weekdays and weekends. Furthermore, data from
December 23rd to December 31st are excluded to prevent distortions
caused by the Christmas holiday season. The evaluations were conducted on a 64-bit Linux server equipped with dual Intel Xeon Gold 6226 CPUs, providing a total of 8 cores running at 2.7 GHz, and 288 GB of RAM. The core components of the column generation algorithm are implemented in C++, with each restricted linear master problem solved using HiGHS version 1.6.0 \citep{huangfu2018parallelizing}. Each instance was subject
to a time limit. If the column generation reaches this limit before
converging, it enters the search phase: the path with the largest
fractional value is fixed at 1, and this process is repeated until a
feasible integer solution is found.

Table~\ref{tab:cg_result_0.5h} summarizes the results across all
instances with a tight 30-minute time limit, that captures the
realities in the field. For context, the third column shows the total
number of trip requests and the number of trips successfully fulfilled
by the current Paratransit system. The current Paratransit system,
which uses a heuristic algorithm provided by RouteMatch given
predetermined shifts, typically requires hours to perform the
scheduling process. On average, this system accommodates 80.75\% of
the total trip requests. 
\revision{In the table, the ``objective" value 
represents the number of requests  served in the final integer 
solution; ``time (s)" denotes the total runtime in seconds for the 
entire algorithm until convergence or the time limit. The 30-minute time limit experiment is
designed to assess the solution quality of the proposed \algo{} when
required to operate within a much shorter time frame. The goal is to
enhance both the planning process and the user experience by
significantly reducing scheduling time compared to the current system.
}

\begin{table}[!ht]
\caption{Performance Comparison Between the Proposed \algo{} and Alternative Approaches With 30-Minute Time Limit.}
    \centering
\begin{tabular}{ccccccccc}
\hline
\multirow{2}{*}{instance} & {number of edges} & total trip requests & \multicolumn{2}{c}{Baseline-CG} & \multicolumn{2}{c}{RF-CG} & \multicolumn{2}{c}{\algo{}} \\ \cline{4-9} 
    & in original graph & / current system & objective & time (s) & objective & time (s) & objective & time (s)\\ \hline
20191202 & 902,025 & 475/365 & 380 & 1,800 & 414 & 1,552 & 437 & 1,534 \\
20191203 & 853,314 & 462/342 & 368 & 1,800 & 407 & 1,800 & 439 & 898  \\
20191204 & 1,258,323 & 561/450 & 427 & 1,800 & 458 & 1,800 & 470 & 1,800 \\
20191205 & 1,060,385 & 515/409 & 387 & 1,800 & 442 & 1,795 & 480 & 1,737 \\
20191206 & 1,035,815 & 509/391 & 401 & 1,800 & 444 & 1,653 & 460 & 1,241 \\
20191209 & 971,703 & 493/427 & 380 & 1,800 & 406 & 1,800 & 461 & 1,594 \\
20191210 & 1,031,748 & 508/419 & 418 & 1,800 & 463 & 1,469 & 478 & 1,571 \\
20191211 & 1,205,055 & 549/443 & 418 & 1,800 & 501 & 1,800 & 450 & 1,800 \\
20191212 & 1,023,638 & 506/421 & 386 & 1,800 & 450 & 1,800 & 435 & 1,800 \\
20191213 & 1,043,973 & 511/398 & 394 & 1,800 & 423 & 1,800 & 479 & 1,231 \\
20191216 & 913,458 & 478/405 & 387 & 1,800 & 431 & 1,123 & 440 & 1,195 \\
20191217 & 831,288 & 456/405 & 436 & 1,793 & 410 & 908  & 435 & 986  \\
20191218 & 898,230 & 474/384 & 374 & 1,800 & 413 & 1,658 & 442 & 1,491 \\
20191219 & 1,072,778 & 518/443 & 402 & 1,800 & 467 & 1,800 & 483 & 1,415 \\
20191220 & 987,539 & 497/363 & 408 & 1,800 & 451 & 1,558 & 456 & 1,397 \\
\hline
average  &  1,005,952 & 501/404 & 398 & 1,800 & 439 & 1,621 & 456 & 1,446\\
\hline
\end{tabular}
    \label{tab:cg_result_0.5h}
\end{table}

Table~\ref{tab:cg_result_0.5h} shows that Baseline-CG hits the
30-minute time limit for almost all instances. This is expected as
Baseline-CG needs much longer time to complete the pricing
subproblems, thus struggling to explore sufficient number of
high-quality paths for the master problem. This results in an average
objective even worse than that of the existing system. In comparison,
both RF-CG and \algo{} operate on reduced graph sizes, this
accelerating the pricing subproblems. They improve solution quality by
10.3\% and 14.6\%, respectively, and reduce runtimes by 9.9\% and
19.7\%, respectively, compared to Baseline-CG. The improvement of
RF-CG over Baseline-CG is in tune with the results reported in
\citet{morabit2023machine}.

Table~\ref{tab:cg_result_0.5h} also indicates that \algo{} provides
significant advantage over RF-CG.  On the one hand, \algo{} achieves a
10\% additional reduction in runtimes and it reaches the time limit on
the three instances fewer than RF-CF. On the other hand, \algo{}
improves the objective by 3.9\% in average compared to RF-CF, and
always improves the solution quality over the current system.  This
contrasts with RF-CG that may fail to improve over the current system
(e.g., instance 20191209). These observations align with the
classification results, as the GNN is able to identify promising edges
more accurately than RF.

To further illustrate of the benefits of \algo{},
Table~\ref{tab:cg_result_1h} presents a comparison of performance
metrics when the time limit is extended to 1 hour. In this case,
neither RF-CG nor \algo{} exceed time limit for any instance, whereas
Baseline-CG continues to reach the time limit in 5 instances.  With
the extended run time, all three CG approaches show an increase in
average objective values compared to the 30-minute time limit as
expected. Baseline-CG, in particular, shows the most significant
improvement, as it is the most affected by the 30-minute time limit
bottleneck. As a result, RF-CG loses its solution quality advantage
over Baseline-CG, now trailing by an average of 3.3\%. In contrast,
\algo{} maintains a 1\% advantage in average objective value while
achieving a remarkable 49.8\% reduction in run time.

\begin{table}[!ht]
\caption{Performance Comparison Between the Proposed \algo{} and Alternative Approaches With 1-Hour Time Limit.}
\centering
\begin{tabular}{ccccccc}
\hline
\multirow{2}{*}{instance} & \multicolumn{2}{c}{Baseline-CG} & \multicolumn{2}{c}{RF-CG} & \multicolumn{2}{c}{\algo{}} \\ 
\cline{2-7} 
     & objective & time (s) & objective & time (s) & objective & time (s)\\ \hline
20191202 & 407 & 3,600 & 414 & 1,552 & 437 & 1,534 \\
20191203 & 450 & 1,999 & 405 & 1,855 & 439 & 898  \\
20191204 & 463 & 3,600 & 517 & 2,288 & 519 & 2,526 \\
20191205 & 493 & 2,576 & 442 & 1,795 & 480 & 1,737 \\
20191206 & 473 & 2,543 & 444 & 1,653 & 460 & 1,241 \\
20191209 & 433 & 3,600 & 431 & 2,133 & 461 & 1,594 \\
20191210 & 479 & 3,335 & 463 & 1,469 & 478 & 1,571 \\
20191211 & 472 & 3,600 & 498 & 1,953 & 508 & 1,871 \\
20191212 & 454 & 3,600 & 459 & 1,900 & 478 & 1,824 \\
20191213 & 484 & 3,000 & 456 & 2,254 & 479 & 1,231 \\
20191216 & 449 & 2,820 & 431 & 1,123 & 440 & 1,195 \\
20191217 & 436 & 1,793 & 410 & 908  & 435 & 986  \\
20191218 & 449 & 3,280 & 413 & 1,658 & 442 & 1,491 \\
20191219 & 499 & 3,170 & 463 & 1,984 & 483 & 1,415 \\
20191220 & 477 & 2,331 & 451 & 1,558 & 456 & 1,397 \\
\hline
average  & 461 & 2,993 & 446 & 1,739 & 466 & 1,501\\
\hline
\end{tabular}
    \label{tab:cg_result_1h}
\end{table}

In conclusion, the combination of improved solution quality and
reduced computational time demonstrates the benefits and robustness of
\algo{}. This is especially valuable for large-scale instances with
high demand variability, making \algo{} highly suitable for real-world
applications where both solution quality and computational efficiency
are critical.

\subsection{Sensitivity Analysis}
\label{section:sensitivity}

This section presents a sensitivity analysis of the impact of
different reduced graph configurations used in GNN training (RGTrn)
and testing (RGTst).  For each training instance, the RGTrn graphs
provide the GNN-based model with labels indicating promising edges for
learning. The quality of these graphs directly impacts the GNN
ability to accurately identify promising edges. During testing, the
RGTst graphs determine the reduced graph on which the column
generation algorithm runs. The size of the RGTst graph represents a
trade-off between solution quality and computation time in the column
generation process.

\begin{figure}[!ht]
    \centering
    \includegraphics[width=\textwidth]{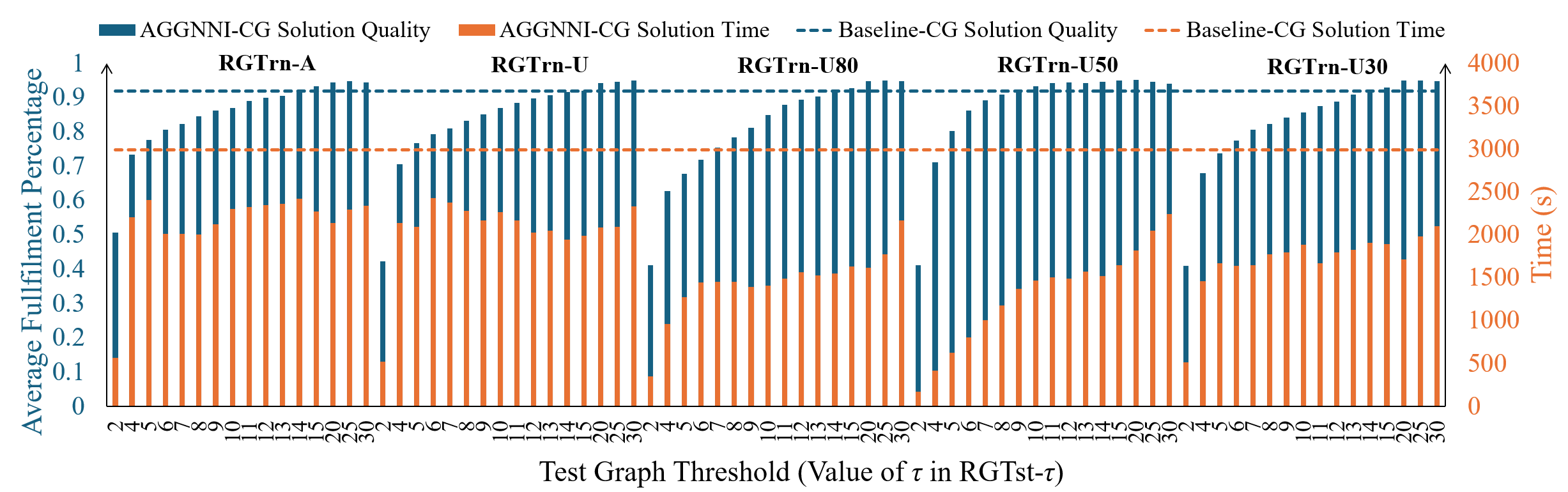}
        \caption{\algo{} Performance Comparison Across Different Types of RGTrn and RGTst Graphs, with Dashed Lines Representing Baseline-CG Metrics and Solid Lines Representing \algo{} Metrics.}
        \label{fig:sensitivity_analysis}
\end{figure}

Figure~\ref{fig:sensitivity_analysis} summarizes the results of the
proposed \algo{} when trained and tested with various types of RGTrn
and RGTst graphs and the 1-hour time limit.  The figure is organized
into five groups, corresponding to the five types of RGTrn
graphs. Within each group, the results for sixteen types of RGTst
graphs are presented. For example, the first bar in
Figure~\ref{fig:sensitivity_analysis} represents the result of using
RGTrn-A as the training graph and RGTst-2 as the test graph. The
height of each bar indicates the average solution quality, while the
height of the orange section reflects the average solution time for
all test instances. The labels of vertical axis on the left and the
right of the figure provide the scale in terms of solution quality and
runtime, respectively.

To ensure a consistent comparison across different test days, which
have varying numbers of trip requests, we calculate the fulfillment
percentage (i.e., the ratio of trips fulfilled to total trip requests)
for each day, and average this across all days to represent solution
quality. Additionally, the average solution quality (blue dashed line)
and solution time (orange dashed line) of Baseline-CG with 1-hour time
limit are included for reference. The following key findings emerge
from the sensitivity analysis in
Figure~\ref{fig:sensitivity_analysis}:

\paragraph{Impact of Training Graphs (RGTrn):}
There is a clear distinction among the five training graph groups. The
average solution times for RGTrn-A and RGTrn-U are higher than those
of the remaining three groups. This suggests that edges selectively
derived from paths with positive usage in solving the RLMP, rather
than from all paths explored during the column generation iterations,
are more suitable to be labeled as promising for supervised learning.
Additionally, \algo{} exhibit high robustness, showing consistent
performance across different RGTrn scales (e.g., RGTrn-U80, RGTrn-U50,
and RGTrn-U30).

\paragraph{Effect of Test Graphs (RGTst):}
In terms of solution quality, smaller RGTst graphs (e.g., RGTst-2 and
RGTst-4) generally result in faster solution times but at the expense
of lower fulfillment percentages.  As the scale of the RGTst graph
increases, the solution quality improves, but this comes with longer
computation times. A more balanced combination of solution quality and
time is observed in mid-range RGTst graphs (e.g., between RGTst-10 and
RGTst-15), where the solution quality approaches or even exceeds that
of Baseline-CG, while the solution time is significantly reduced.

\paragraph{Comparison with Baseline-CG Performance:}
Across all RGTrn types, the proposed \algo{} consistently achieves a
higher average fulfillment percentage than Baseline-CG, while
drastically reducing the solution time. This demonstrates the
efficiency of the GNN-based graph reduction approach.  Particularly
for RGTrn-U80, RGTrn-U50 and RGTrn-U30, the combination of solution
quality and solution time shows a substantial improvement over
Baseline-CG, making these configurations highly suitable for practical
applications.

In summary, the sensitivity analysis confirms that the proposed
\algo{} is robust to variations in the scale of training graphs and
achieves a desirable balance between solution quality and time for
mid-scale test graphs. These results suggest that \algo{} is
well-suited for large-scale, real-world \pb{} problems, where both
computational efficiency and solution accuracy are critical.

\section{Conclusions}
\label{sec:conclusion}

This paper introduced the Joint Rider Trip Planning and Crew Shift
Scheduling Problem (\pb{}) to avoid separating two important
components of various mobility systems. Moreover, to meet the
computational challenges of \pb{}, the paper proposed a novel solution
method, called Attention and Gated GNN-Informed Column Generation
(\algo{}) that hybridizes column generation and machine learning.
The key idea underlying \algo{} is to reduce the size of the graphs
explored by the pricing subproblems by only exploring promising edges,
accelerating the most time-consuming component of the column
generation.  These edges are obtained by supervised learning, using a
novel graph neural network that features both a multi-head attention
mechanism and a gated architecture. The proposed GNN is ideally suited
to cater for the different input sizes encountered during daily
operations, where the number of requests, their locations, and request
times vary from day to day.

The benefits of \algo{} was demonstrated on a real-world paratransit
systems, where it was able to condense the graphs for the pricing
subproblems by an order of magnitude. This led to a significant
acceleration of the conventional column generation process. For
complex instances where traditional column generation methods
falter, \algo{} successfully delivers high-quality solutions within
a practical timeframe. The methodology underlying \algo{} is also
general and should apply to other applications as well.

There are several avenues for expanding upon this research. First,
while this study leverages a reservation-based on-demand travel system
for performance validation, it would be highly interesting to apply a
similar approach to real-world multimodal transit systems. Second, it
would be interesting to determine whether it is possible to generalize
\algo{} to handle seamlessly both weekdays, weekends, and holidays,
simplifying the deployment in practice. \revision{Third, extending the framework to support online rescheduling would address operational stochasticity in transportation systems. This would require incorporating real-time data streams to dynamically adjust schedules in response to disruptions like traffic congestion, vehicle breakdowns, and passenger no-shows.}

\section*{Declaration of generative AI and AI-assisted technologies in the writing process}
During the preparation of this work the authors used ChatGPT in order to
improve readability. After using this service, the authors reviewed and
edited the content as needed and take full responsibility for the content
of the publication.



\newpage
\appendix 

\section{Column Generation Details}
\label{app:cg_details}
\revision{For the completeness of this paper, this section provides the pseudocode of the column generation algorithm proposed in Section~\ref{sec:column_generation}.}

\SetKwFunction{ColGenFixed}{ColumnGenerationFixed}
\SetKwFunction{SolvePricing}{SolvePricing}
\SetKwFunction{LabelExt}{LabelExtension}
\SetKwFunction{RLMP}{SolveRLMP}
\SetKwFunction{FindMaxFracCol}{FindMaxFractionalColumn}
\SetKwFunction{FixColToOne}{FixColumnToOne}
\SetKwFunction{BackToDepot}{BackToDepot}
\SetKwFunction{DropOff}{DropOff}
\SetKwFunction{PickUp}{PickUp}
\SetKwFunction{TimeUpdate}{TimeUpdate}
\SetKwFunction{TimeFeasible}{TimeFeasible}
\SetKwFunction{ReachableUpdate}{ReachableUpdate}
\SetKwFunction{BacktraceRoute}{BacktraceRoute}
\SetKwFunction{CheckDominance}{CheckDominance}
\SetKwFunction{ExtractMin}{ExtractMin}
\SetKwFunction{Insert}{Insert}
\SetKwProg{Proc}{Procedure}{}{}

\subsection{Column Generation Outer Loop}

\revision{Algorithm~\ref{alg:outer} outlines a branching heuristic to solve an integer programming formulation for \pb{}. The goal is to iteratively construct an integral feasible solution $\lambda$ by solving a series of restricted master problems with dynamically added columns and branching on fractional solutions.}
{\color{blue}
\begin{algorithm}[!ht]
    \caption{Column Generation Outer Loop (Branching Heuristic in Section~\ref{sec:feasible_solution})}
    \label{alg:outer}
    \SetKwInOut{Input}{Input}
    \SetKwInOut{Output}{Output}
    \SetKw{Break}{break}
    \SetKw{Continue}{continue}
    \Input{Set of trip requests $R$, Fleet $F$, Set of driver shift candidates $\Phi$,\\
        Graph $g=(N^g,E^g)$ with travel time $t_{ij}$ associated with each edge $(i, j) \in E^g$}
    \Output{An integer feasible solution $\lambda$}

    $\Omega' \gets \emptyset$ \; 
    $\lambda \gets \text{None}$ \;
    $\mathit{fixed\_columns} \gets \emptyset$ \;

    \While{$\lambda$ \text{is None}}{
        $\lambda' \gets \ColGenFixed(\Omega', R, F, \Phi, \mathit{fixed\_columns}, g)$\; \label{line:call_cg_fixed}
        \If{$\lambda'$ \text{is integer}}{
            $\lambda \gets \lambda'$\;
            \Break\;
        }

        $\theta_{\mathit{fix}} \gets \FindMaxFracCol(\lambda')$\; \label{line:find_frac}
        $\mathit{fixed\_columns} \gets \mathit{fixed\_columns} \cup \{\theta_{\mathit{fix}}\}$\; \label{line:fix_col}
    }
    \Return{$\lambda$}\;
\end{algorithm}
}

\revision{Lines 1--3 initialize the algorithm. The fixed columns set is used to record branching decisions, where specific columns are constrained to be 1 in subsequent iterations. Lines 4--10 form the core iterative loop. In each iteration, the algorithm finds new columns through the \texttt{ColumnGenerationFixed} procedure using the current problem state and fixed columns. If the resulting solution $\lambda'$ is integer, it is assigned to $\lambda$ and the algorithm terminates. Otherwise, the largest fractional column $\theta_{\text{fix}}$ is identified using the \texttt{FindMaxFractionalColumn} routine, and this column is added to the set of fixed columns to guide the next iteration.
Line 11 returns the final integer feasible solution $\lambda$ once the loop has converged.}

\subsection{Column Generation with Fixed Columns}

\revision{Algorithm~\ref{alg:cg_fixed} describes the inner loop of the column generation procedure, which solves a restricted linear master problem (RLMP) iteratively while incorporating a fixed set of columns due to branching constraints. The algorithm starts with an initial set of columns $\Omega'$ and repeatedly solves the RLMP using these columns and the current set of fixed variables. The RLMP solution returns dual variables associated with the trip constraints and resource constraints.}

\begin{algorithm}[!ht]
    \caption{Column Generation with Fixed Columns [\texttt{ColumnGenerationFixed}]}
    \label{alg:cg_fixed}
    \SetKwInOut{Input}{Input}
    \SetKwInOut{Output}{Output}

    \Input{Current set of columns $\Omega'$, Set of trip requests $R$, Fleet $F$, Set of driver shift candidates $\Phi$, \\
    Set of fixed columns $\mathit{fixed\_columns}$, Graph $g=(N^g,E^g)$}
    \Output{RLMP solution $\lambda'$}

    \While{True}{
        $\lambda', \{\pi_r'\}_{r \in R}, \sigma' \gets \RLMP(\Omega', \mathit{fixed\_columns}, R, F)$ \tcp*{Solve the RLMP (Figure~\ref{fig:cg_mp}) with columns $\Omega'$ and additional constraints on $\mathit{fixed\_columns}$}\label{line:solve_rlmp} 
        $\mathcal{C}_{new} \gets \emptyset$\;
        \For{$\phi = (dr_s, dr_e) \in \Phi$}{ \label{line:loop_shifts}
             $\mathcal{C}_{\phi} \gets \SolvePricing(\{\pi_r'\}_{r \in R}, \sigma', \phi, g)$\; \label{line:call_pricing}
             $\mathcal{C}_{\mathit{new}} \gets \mathcal{C}_{\mathit{new}} \cup \mathcal{C}_{\phi}$\;
        }
        \If{$\mathcal{C}_{\mathit{new}} = \emptyset$}{ \label{line:check_new_cols}
            \Return{$\lambda'$} \tcp*{Terminate when the Pricing Subproblem yields no new columns}
        }
        $\Omega' \gets \Omega' \cup \mathcal{C}_{\mathit{new}}$\; \label{line:add_cols}
   }
\end{algorithm}

\revision{
Lines 1--2 initialize the iterative process. In each iteration, the RLMP is solved with the current column pool $\Omega'$ and the fixed column constraints, yielding a primal solution $\lambda'$ and associated dual variables. Lines 4--6 loop over all shift candidates in $\Phi$ and solve the pricing subproblem for each shift $\phi$, using the current dual variables to identify columns with negative reduced cost. Any generated columns are added to the set $C_{\text{new}}$.

If no new columns are generated in an iteration, as checked in line 7, the algorithm terminates and returns the current RLMP solution $\lambda'$. Otherwise, line 9 updates the column set $\Omega'$ by including the new columns, and the process repeats until optimality or time limit is reached.}

\subsection{Solving the Pricing Subproblem}

\revision{Algorithm~\ref{alg:pricing} solves the pricing subproblem for a given driver shift $\phi$ by identifying columns with negative reduced cost. It employs a dynamic programming algorithm on a graph (an original or a reduced graph in this study), where each label represents a partial route that originates at the depot and terminates at a specific node. A label is defined as $L = (i, c, O, Re, S, t_a, t_s, t_d, \mathit{prev})$, where $i$ denotes the current node, $c$ is the reduced cost of the partial path, $O$ is the set of open trip requests (i.e., requests that have been picked up but not yet dropped off), $Re$ is the set of reachable trip requests, and $S$ is the set of already served requests. The variables $t_a$, $t_s$, and $t_d$ represent the arrival, service start, and departure times at node $i$, respectively. \revisiontwo{These label components collectively track the consumption of resources (time, capacity) and ensure that constraints \eqref{eq:same_vehicle} - \eqref{eq:arc_model_z} are satisfied during path construction.} The field $\mathit{prev}$ stores a reference to the preceding label in the route, enabling backtracking to reconstruct complete paths.}

\begin{algorithm}[!ht]
    \caption{Solve the Pricing Subproblem (Section \ref{sec:pricing_subproblem}) [\texttt{SolvePricing}]}
    \label{alg:pricing}
    \SetKwInOut{Input}{Input}
    \SetKwInOut{Output}{Output}
    \SetKw{Break}{break}
    \SetKw{Continue}{continue}
    \Input{Dual variables $\{\pi_r\}_{r \in R}$ (requests) and $\sigma$ (fleet size), Driver shift $\phi=(dr_s, dr_e)$, Graph $g=(N^g,E^g)$}
    \Output{Set of new columns $\mathcal{C}_{\phi}$ with negative reduced cost for shift $\phi$}

    $\mathit{L_{pool}} \gets \text{Empty priority queue ordered by reduced cost}$\;
    $\mathit{Processed\_Labels}_i \gets \emptyset$ for every node $i \in N^g$\;
    $\mathit{Columns\_Found} \gets \emptyset$\;
    $\mathit{max\_columns} \gets \text{predefined limit}$\;
    $\mathit{max\_labels\_per\_node} \gets \text{predefined limit}$\;

    $\mathit{L_0} \gets (i=0, c=-\sigma, O=\emptyset, Re=P, S=\emptyset, t_a=dr_s, t_s=dr_s, t_d=dr_s, \mathit{prev}=\emptyset)$\; \label{line:init_label}
    $\mathit{L_{pool}}.\Insert(\mathit{L_0})$\;
    $\mathit{Processed\_Labels}_0 \gets \mathit{Processed\_Labels}_0 \cup \{\mathit{L_0}\}$\;

    \While{$(\mathit{L_{pool}} \text{ is not empty}) \land (|\mathit{Columns\_Found}| < \mathit{max\_columns})$}{
        $L = (i, c, O, Re, S, t_a, t_s, t_d, \mathit{prev}) \gets \mathit{L_{pool}}.\ExtractMin()$ \tcp*{label with the lowest cost}
        $\mathit{Generated\_Labels} \gets \LabelExt(L, \{\pi_r\}_{r \in R}, \phi, g)$\; \label{line:call_labelext}

        \For{$L' = (j, c', O', Re', S', t_a', t_s', t_d', \mathit{prev}') \in \mathit{Generated\_Labels}$}{
            \If{$j = 2n+1$}{\label{line:check_depot}
                 $\theta \gets \BacktraceRoute(L')$ \tcp*{Backtrace route upon reaching destination depot}
                 $\mathit{Columns\_Found} \gets \mathit{Columns\_Found} \cup \{\theta\}$\;
                 \If{$|\mathit{Columns\_Found}| \ge \mathit{max\_columns}$}{\Break\;}
            }
            \Else{
                 $\mathit{dominated} \gets \text{False}$\;
                 $\mathit{Labels\_to\_Remove} \gets \emptyset$\;
                 \For{$L'' \in \mathit{Processed\_Labels}_j$}{
                    \If{\CheckDominance{$L', L''$}}{
                        $\mathit{dominated} \gets \text{True};$ \text{break;}
                    }
                    \If{\CheckDominance{$L'', L'$}}{
                        $\mathit{Labels\_to\_Remove} \gets \mathit{Labels\_to\_Remove} \cup \{L''\}$ \tcp*{Remove $L''$ if dominated by $L'$}
                    }
                 }
                 \If{$\mathit{dominated} = \text{True}$}{\Continue \tcp*{ Skip if $L'$ is dominated by any $L'' \in \mathit{Processed\_Labels}_j$}}
                 
                 $\mathit{Processed\_Labels}_j \gets \mathit{Processed\_Labels}_j \setminus \mathit{Labels\_to\_Remove}$\;
                 $\mathit{L_{pool}}.\Insert(L')$\;
                 $\mathit{Processed\_Labels}_j \gets \mathit{Processed\_Labels}_j \cup \{L'\}$\;

                 \If{$|\mathit{Processed\_Labels}_j| > \mathit{max\_labels\_per\_node}$}{
                      Remove label(s) with the highest cost from $\mathit{Processed\_Labels}_j$\;
                 }
            }
        }
    }
    \Return{$\mathit{Columns\_Found}$}\;

    \Proc{\CheckDominance{$L', L''$}}{
        Let $L' = (j, c', O', Re', S', t_a', t_s', t_d')$\;
        Let $L'' = (j, c'', O'', Re'', S'', t_a'', t_s'', t_d'')$\;
        $\mathit{dominated} \gets (O'' \subseteq O') \land (t_d'' \le t_d') \land (c'' \le c') \land (|S''| \ge |S'|)$ \tcp*{$L'$ dominated by $L''$}
        \Return{$\mathit{dominated}$}\;
    }
\end{algorithm}

\revision{
Lines 1--4 initialize key data structures, including the priority queue $L_{\text{pool}}$, a processed label pool for each node, and the output set of columns $\mathit{Columns\_Found}$. Line 5 creates the initial label at the start depot and inserts it into the queue. Lines 6–8 define the main loop, which continues until the queue is empty or a predefined number of columns have been found. At each iteration, the label with the lowest cost is extracted and extended via the \texttt{LabelExtension} function (line 9). If the label reaches the destination depot (line 10), a feasible route is reconstructed using \texttt{BacktraceRoute}, and the resulting column is added to the solution set (lines 11–13). The loop terminates early if the maximum number of columns is reached.

If the current label does not reach the destination depot, lines 15–28 perform dominance checks to reduce the number of stored labels. A label $L'$ is compared against previously processed labels at the same node to detect dominance relationships. If $L'$ is dominated, it is skipped; if it dominates others, those are marked for removal. The surviving label is added to the processed set, and the label pool is updated. If the number of labels at a node exceeds a predefined limit, the labels with the highest cost are pruned.

Lines 29–31 define the dominance-checking procedure used above. A label $L'$ is considered dominated by another label $L''$ if $L''$ arrives earlier, has lower cost, and represents a superset of visited requests. 
}

\subsection{Label Extension Details}

\revision{
Algorithm~\ref{alg:label_ext} implements the label extension procedure used within the pricing subproblem. Given a current label $L = (i, c, O, Re, S, t_a, t_s, t_d, \mathit{prev})$, the algorithm generates feasible extensions by considering depot returns, drop-off actions, and pickup opportunities for unserved trip requests. Each extension results in a new label, which is added to the set of generated labels if it satisfies temporal and feasibility constraints.
}

\begin{algorithm}[!ht]
    \caption{Label Extension [\texttt{LabelExtension}]}
    \label{alg:label_ext}
    \SetKwInOut{Input}{Input}
    \SetKwInOut{Output}{Output}
    \SetKw{Break}{break}
    \SetKw{Continue}{continue}
     \Input{Current Label $L=(i, c, O, Re, S, t_a, t_s, t_d, \mathit{prev})$, Dual variables $\{\pi_r\}_{r \in R}$, Driver shift $\phi=(dr_s, dr_e)$, Graph $g=(N^g,E^g)$}
    \Output{Set of extended labels $\mathit{Generated\_Labels}$}

    $\mathit{Generated\_Labels} \gets \emptyset$\;

    \If{$(i \neq 0) \land (|O| = 0) \land (c < 0)$}{
        $j \gets 2n+1$ \tcp*{Back to depot}
        $(t_a', t_s', t_d') \gets \TimeUpdate(t_d, i, j)$\;
        \If{\TimeFeasible{$t_a', j, dr_e$}}{
            $c' \gets c$; $O' \gets \emptyset$; $Re' \gets \emptyset$; $S' \gets S$\;
            $L' \gets (j, c', O', Re', S', t_a', t_s', t_d', \mathit{prev'}=L)$\;
            $\mathit{Generated\_Labels} \gets \mathit{Generated\_Labels} \cup \{L'\}$\;
        }
    }

    \For{$r \in O$}{
         $j \gets n + r$ \tcp*{Drop-off node $n + r$ for request $r$}
         \If{$(i,j) \notin E^g$}{\Continue}
         $(t_a', t_s', t_d') \gets \TimeUpdate(t_d, i, j)$\;
         \If{\TimeFeasible{$t_a', j, dr_e$}}{
             $c' \gets c$;
             $O' \gets O \setminus \{r\}$;
             $S' \gets S$;
             $Re' \gets \ReachableUpdate(Re, j)$\;
             $L' \gets (j, c', O', Re', S', t_a', t_s', t_d', \mathit{prev'}=L)$\;
             $\mathit{Generated\_Labels}  \gets \mathit{Generated\_Labels} \cup \{L'\}$;
         }
    }

    \For{$r \in Re$}{
          $j \gets r$ \tcp*{Pickup node $r$ for request $r$}
          \If{$(i,j) \notin E^g$}{\Continue}
          \If{$(|O| < C) \land (r \notin O) \land (r \notin S)$}{
             $(t_a', t_s', t_d') \gets \TimeUpdate(t_d, i, j)$\;
             \If{\TimeFeasible{$t_a', j, dr_e$}}{
                 $c' \gets c - \pi_r$;
                 $O' \gets O \cup \{r\}$;
                 $S' \gets S \cup \{r\}$;
                 $Re' \gets \ReachableUpdate(Re, j)$\;
                 $L' \gets (j, c', O', Re', S', t_a', t_s', t_d', \mathit{prev'}=L)$\;
                $\mathit{Generated\_Labels}  \gets \mathit{Generated\_Labels} \cup \{L'\}$;             }
         }
    }

    \Return{$\mathit{Generated\_Labels}$}\;

    \Proc{\TimeUpdate{$t_d^{\mathit{prev}}, i, j$}}{
        $t_a' = t_d^{\mathit{prev}} + t_{ij}$;
        $t_s' = \max \{t_a', a_j\}$;
        $t_d' = t_s' + s_j$\;
        \Return{$(t_a', t_s', t_d')$}\;
    }

     \Proc{\ReachableUpdate{$Re^{\mathit{prev}}, \mathit{current\_node}$}}{
          $Re' \gets \{r \in Re^{\mathit{prev}} \mid \TimeFeasible(\mathit{current\_node}, r) \land \TimeFeasible(r, r+n) \land \TimeFeasible(r+n, 2n+1)\}$\;
          \Return{$Re'$}\;
     }
    \Proc{\TimeFeasible{$t_a', j, dr_e$}}{
         \Return{$(t_a' \le b_j) \land (t_a' \le dr_e)$}\; 
    }
\end{algorithm}

\revision{
Lines 1–8 handle the case where the route may return to the depot. If the label is not at the depot and has no open requests, and the reduced cost is negative, a depot return is considered. The algorithm computes updated timing information and checks time feasibility before appending the corresponding extended label.

Lines 9–17 iterate over all served requests to evaluate drop-off possibilities. For each request $r$, the corresponding drop-off node $j = n + r$ is checked for edge existence and temporal feasibility. If feasible, a new label is generated with updated cost, state sets, and reachable request sets.

Lines 18–27 focus on pickup actions for all reachable trip requests. For each candidate $r \in Re$, the pickup node is evaluated under precedence and capacity constraints (e.g., the number of open requests must be below vehicle capacity, and a request cannot be picked up again if already served). If all feasibility conditions are satisfied, a new label is constructed with the updated attributes and added to the output.

The subprocedures at the end define utility functions. \texttt{TimeUpdate} (lines 29–31) computes updated arrival, service, and departure times. \texttt{ReachableUpdate} (lines 32–34) filters reachable requests based on temporal feasibility from the new node. Finally, \texttt{TimeFeasible} (lines 35–36) checks whether a proposed arrival time is within the time window of the target node.
}

\revisiontwo{
\section{Feature Importance Analysis}\label{sec:feature_importance}

To evaluate the relative importance of features used in the GNN model, this paper extends existing feature importance methods for GNNs. To the author's knowledge, current methods like GNNExplainer \citep{ying2019gnnexplainer} only support importance ranking for node features. To address this limitation, this paper enhances the GNNExplainer framework to incorporate edge feature masking alongside node feature masking.

Table \ref{tab:gnnexplainer} introduces the notations used in the enhanced GNNExplainer. The following augmented loss function that integrates edge feature masking is proposed:
\begin{align}
\mathcal{L}_{\text{enhanced}} = \mathcal{L}_{\text{pred}}(\hat{y}, y) + \mathcal{L}_{\text{node}} + \mathcal{L}_{\text{edge}},
\end{align}
where
\begin{align}
\mathcal{L}_{\text{edge}} = & \lambda_{\text{edge\_size}} \cdot \frac{1}{|B_{\text{edge}}|} \sum_{i \in B_{\text{edge}}} \sigma(M_{\text{edge}}[i]) \nonumber\\ 
&+ \lambda_{\text{edge\_ent}} \cdot \frac{1}{|B_{\text{edge}}|} \sum_{i \in B_{\text{edge}}} H(\sigma(M_{\text{edge}}[i])),\\
\mathcal{L}_{\text{node}} = & \lambda_{\text{node\_size}} \cdot \frac{1}{|B_{\text{node}}|} \sum_{i \in B_{\text{node}}} \sigma(M_{\text{node}}[i]) \nonumber \\
& \lambda_{\text{node\_ent}} \cdot \frac{1}{|B_{\text{node}}|} \sum_{i \in B_{\text{node}}} H(\sigma(M_{\text{node}}[i]), \\
\mathcal{L}_{\text{pred}} =  &
-\sum_i [y_i \log(\hat{y}_i) + (1-y_i)\log(1-\hat{y}_i)],\\
H(p) = &  -p\log(p+\epsilon) - (1-p)\log(1-p+\epsilon).
\end{align}

\begin{table}[!ht]
\centering
\caption{Enhanced GNNExplainer loss function parameters}
\begin{tabular}{llc}
\hline
\textbf{Symbol} & \textbf{Definition} & \textbf{Value} \\
\hline
$M_{\text{node}}$ & Learnable node feature importance parameters & - \\
$M_{\text{edge}}$ & Learnable edge feature importance parameters & - \\
$B_{\text{node}}$ & Boolean mask indicating node features with non-zero gradients & - \\
$B_{\text{edge}}$ & Boolean mask indicating edge features with non-zero gradients & - \\
$\sigma(\cdot)$ & Sigmoid function & - \\
$\lambda_{\text{node\_size}}$ & Node sparsity coefficient & 1.0 \\
$\lambda_{\text{node\_ent}}$ & Node entropy coefficient & 0.1 \\
$\lambda_{\text{edge\_size}}$ & Edge sparsity coefficient & 1.0 \\
$\lambda_{\text{edge\_ent}}$ & Edge entropy coefficient & 0.1 \\
$\epsilon$ & Numerical stability constant & $10^{-15}$ \\
\hline
\end{tabular}
\label{tab:gnnexplainer}
\end{table}

\subsection{Feature Importance Results and Implications} \label{}

To determine the final importance score for each feature, the learned mask values across all nodes/edges were averaged and then summed across the 15 test instances. Figure \ref{fig:feat-imp} presents the ranking of feature importance in the GNN model.

\begin{figure}[!ht]
\centering
\includegraphics[width=0.8\linewidth]{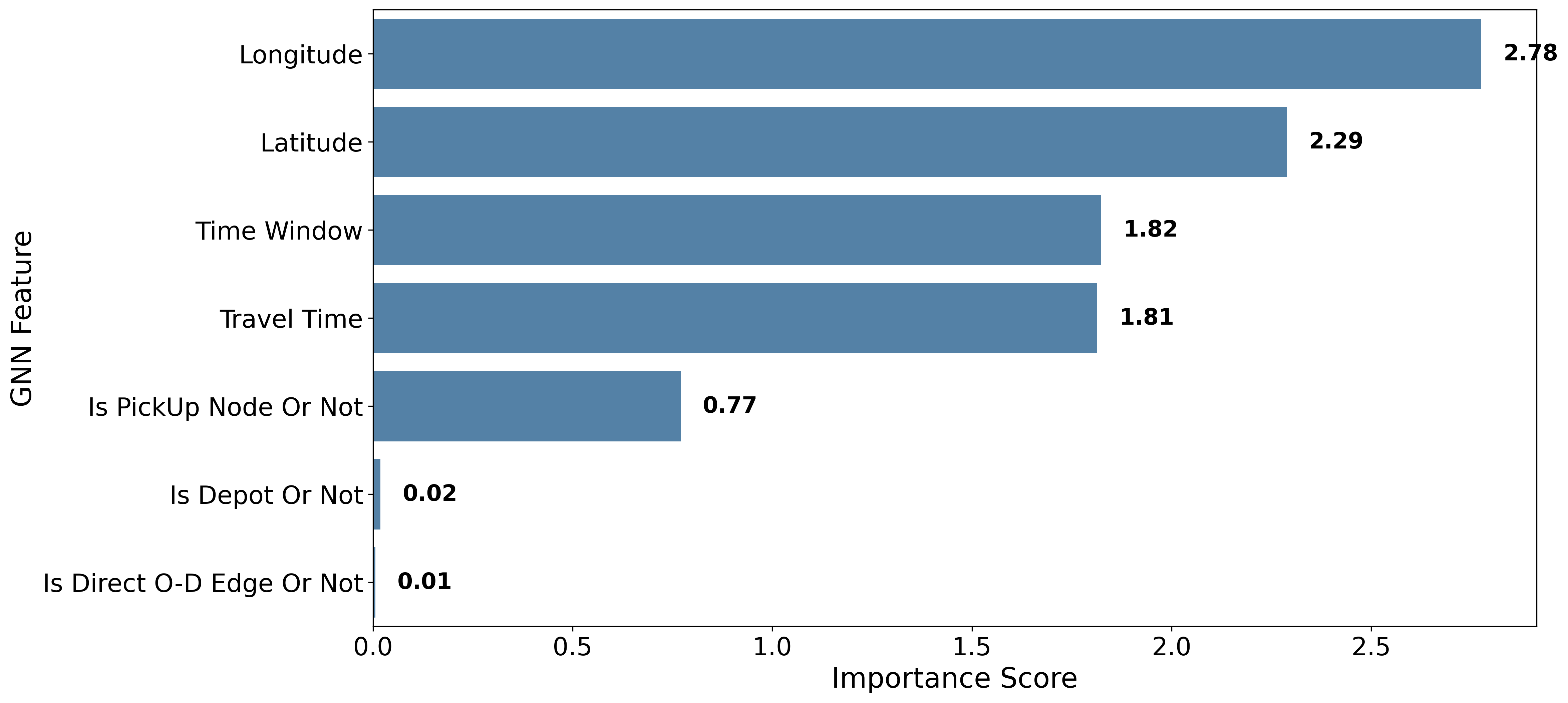}
\caption{Feature Importance on the GNN Model}
\label{fig:feat-imp}
\end{figure}

The analysis reveals that spatial information, specifically the geographic coordinates, contributes most significantly to the model's performance. Longitude (2.78) and latitude (2.29) attained the highest importance scores, highlighting the fundamental role of spatial positioning in transportation routing decisions. Time window constraints (1.82) and travel time (1.81) follow closely with nearly identical importance scores, demonstrating the complementary influence of temporal features on the GNN model.
}

\clearpage
\bibliographystyle{tre}
\bibliography{reference}

\begin{thebibliography}{50}
\expandafter\ifx\csname natexlab\endcsname\relax\def\natexlab#1{#1}\fi
\providecommand{\url}[1]{\texttt{#1}}
\providecommand{\href}[2]{#2}
\providecommand{\path}[1]{#1}
\providecommand{\DOIprefix}{doi:}
\providecommand{\ArXivprefix}{arXiv:}
\providecommand{\URLprefix}{URL: }
\providecommand{\Pubmedprefix}{pmid:}
\providecommand{\doi}[1]{\href{http://dx.doi.org/#1}{\path{#1}}}
\providecommand{\Pubmed}[1]{\href{pmid:#1}{\path{#1}}}
\providecommand{\bibinfo}[2]{#2}
\ifx\xfnm\relax \def\xfnm[#1]{\unskip,\space#1}\fi
\bibitem[{Akiba et~al.(2019)Akiba, Sano, Yanase, Ohta and Koyama}]{akiba2019optuna}
\bibinfo{author}{Akiba, T.}, \bibinfo{author}{Sano, S.}, \bibinfo{author}{Yanase, T.}, \bibinfo{author}{Ohta, T.}, \bibinfo{author}{Koyama, M.}, \bibinfo{year}{2019}.
\newblock \bibinfo{title}{Optuna: A next-generation hyperparameter optimization framework}, in: \bibinfo{booktitle}{Proceedings of the 25th ACM SIGKDD international conference on knowledge discovery \& data mining},  \bibinfo{pages}{2623--2631}.
\bibitem[{Amberg et~al.(2019)Amberg, Amberg and Kliewer}]{amberg2019robust}
\bibinfo{author}{Amberg, B.}, \bibinfo{author}{Amberg, B.}, \bibinfo{author}{Kliewer, N.}, \bibinfo{year}{2019}.
\newblock \bibinfo{title}{Robust efficiency in urban public transportation: Minimizing delay propagation in cost-efficient bus and driver schedules}.
\newblock \bibinfo{journal}{Transportation Science} \bibinfo{volume}{53}, \bibinfo{pages}{89--112}.
\bibitem[{Andrade-Michel et~al.(2021)Andrade-Michel, R{\'\i}os-Sol{\'\i}s and Boyer}]{andrade2021vehicle}
\bibinfo{author}{Andrade-Michel, A.}, \bibinfo{author}{R{\'\i}os-Sol{\'\i}s, Y.A.}, \bibinfo{author}{Boyer, V.}, \bibinfo{year}{2021}.
\newblock \bibinfo{title}{Vehicle and reliable driver scheduling for public bus transportation systems}.
\newblock \bibinfo{journal}{Transportation Research Part B: Methodological} \bibinfo{volume}{145}, \bibinfo{pages}{290--301}.
\bibitem[{Bengio et~al.(2021)Bengio, Lodi and Prouvost}]{bengio2021machine}
\bibinfo{author}{Bengio, Y.}, \bibinfo{author}{Lodi, A.}, \bibinfo{author}{Prouvost, A.}, \bibinfo{year}{2021}.
\newblock \bibinfo{title}{Machine learning for combinatorial optimization: a methodological tour d’horizon}.
\newblock \bibinfo{journal}{European Journal of Operational Research} \bibinfo{volume}{290}, \bibinfo{pages}{405--421}.
\bibitem[{Bornd{\"o}rfer et~al.(2017)Bornd{\"o}rfer, Schulz, Seidl and Weider}]{borndorfer2017integration}
\bibinfo{author}{Bornd{\"o}rfer, R.}, \bibinfo{author}{Schulz, C.}, \bibinfo{author}{Seidl, S.}, \bibinfo{author}{Weider, S.}, \bibinfo{year}{2017}.
\newblock \bibinfo{title}{Integration of duty scheduling and rostering to increase driver satisfaction}.
\newblock \bibinfo{journal}{Public Transport} \bibinfo{volume}{9}, \bibinfo{pages}{177--191}.
\bibitem[{Braekers et~al.(2014)Braekers, Caris and Janssens}]{braekers2014exact}
\bibinfo{author}{Braekers, K.}, \bibinfo{author}{Caris, A.}, \bibinfo{author}{Janssens, G.K.}, \bibinfo{year}{2014}.
\newblock \bibinfo{title}{Exact and meta-heuristic approach for a general heterogeneous dial-a-ride problem with multiple depots}.
\newblock \bibinfo{journal}{Transportation Research Part B: Methodological} \bibinfo{volume}{67}, \bibinfo{pages}{166--186}.
\bibitem[{Bresson and Laurent(2017)}]{bresson2017residual}
\bibinfo{author}{Bresson, X.}, \bibinfo{author}{Laurent, T.}, \bibinfo{year}{2017}.
\newblock \bibinfo{title}{Residual gated graph convnets}.
\newblock \bibinfo{journal}{arXiv preprint arXiv:1711.07553}.
\bibitem[{Brody et~al.(2021)Brody, Alon and Yahav}]{brody2021attentive}
\bibinfo{author}{Brody, S.}, \bibinfo{author}{Alon, U.}, \bibinfo{author}{Yahav, E.}, \bibinfo{year}{2021}.
\newblock \bibinfo{title}{How attentive are graph attention networks?}
\newblock \bibinfo{journal}{arXiv preprint arXiv:2105.14491}.
\bibitem[{Cordeau(2006)}]{cordeau2006branch}
\bibinfo{author}{Cordeau, J.F.}, \bibinfo{year}{2006}.
\newblock \bibinfo{title}{A branch-and-cut algorithm for the dial-a-ride problem}.
\newblock \bibinfo{journal}{Operations Research} \bibinfo{volume}{54}, \bibinfo{pages}{573--586}.
\bibitem[{Cordeau and Laporte(2003)}]{cordeau2003tabu}
\bibinfo{author}{Cordeau, J.F.}, \bibinfo{author}{Laporte, G.}, \bibinfo{year}{2003}.
\newblock \bibinfo{title}{A tabu search heuristic for the static multi-vehicle dial-a-ride problem}.
\newblock \bibinfo{journal}{Transportation Research Part B: Methodological} \bibinfo{volume}{37}, \bibinfo{pages}{579--594}.
\bibitem[{Cordeau and Laporte(2007)}]{cordeau2007dial}
\bibinfo{author}{Cordeau, J.F.}, \bibinfo{author}{Laporte, G.}, \bibinfo{year}{2007}.
\newblock \bibinfo{title}{The dial-a-ride problem: models and algorithms}.
\newblock \bibinfo{journal}{Annals of operations research} \bibinfo{volume}{153}, \bibinfo{pages}{29--46}.
\bibitem[{{Cordeau, Jean-Francois and Groupe d'{\'e}tudes et de recherche en analyse des d{\'e}cisions (Montr{\'e}al, Qu{\'e}bec)}(2000)}]{cordeau2000vrp}
\bibinfo{author}{{Cordeau, Jean-Francois and Groupe d'{\'e}tudes et de recherche en analyse des d{\'e}cisions (Montr{\'e}al, Qu{\'e}bec)}}, \bibinfo{year}{2000}.
\newblock \bibinfo{title}{The VRP with time windows}.
\newblock \bibinfo{publisher}{Citeseer}.
\bibitem[{Cubillos et~al.(2009)Cubillos, Urra and Rodr{\'\i}guez}]{cubillos2009application}
\bibinfo{author}{Cubillos, C.}, \bibinfo{author}{Urra, E.}, \bibinfo{author}{Rodr{\'\i}guez, N.}, \bibinfo{year}{2009}.
\newblock \bibinfo{title}{Application of genetic algorithms for the darptw problem}.
\newblock \bibinfo{journal}{International Journal of Computers Communications \& Control} \bibinfo{volume}{4}, \bibinfo{pages}{127--136}.
\bibitem[{Feng et~al.(2024)Feng, Lusby, Zhang, Tao, Zhang and Peng}]{feng2024branch}
\bibinfo{author}{Feng, T.}, \bibinfo{author}{Lusby, R.M.}, \bibinfo{author}{Zhang, Y.}, \bibinfo{author}{Tao, S.}, \bibinfo{author}{Zhang, B.}, \bibinfo{author}{Peng, Q.}, \bibinfo{year}{2024}.
\newblock \bibinfo{title}{A branch-and-price algorithm for integrating urban rail crew scheduling and rostering problems}.
\newblock \bibinfo{journal}{Transportation Research Part B: Methodological} \bibinfo{volume}{183}, \bibinfo{pages}{102941}.
\bibitem[{Garaix et~al.(2010)Garaix, Artigues, Feillet and Josselin}]{garaix2010vehicle}
\bibinfo{author}{Garaix, T.}, \bibinfo{author}{Artigues, C.}, \bibinfo{author}{Feillet, D.}, \bibinfo{author}{Josselin, D.}, \bibinfo{year}{2010}.
\newblock \bibinfo{title}{Vehicle routing problems with alternative paths: An application to on-demand transportation}.
\newblock \bibinfo{journal}{European Journal of Operational Research} \bibinfo{volume}{204}, \bibinfo{pages}{62--75}.
\bibitem[{Garaix et~al.(2011)Garaix, Artigues, Feillet and Josselin}]{garaix2011optimization}
\bibinfo{author}{Garaix, T.}, \bibinfo{author}{Artigues, C.}, \bibinfo{author}{Feillet, D.}, \bibinfo{author}{Josselin, D.}, \bibinfo{year}{2011}.
\newblock \bibinfo{title}{Optimization of occupancy rate in dial-a-ride problems via linear fractional column generation}.
\newblock \bibinfo{journal}{Computers \& Operations Research} \bibinfo{volume}{38}, \bibinfo{pages}{1435--1442}.
\bibitem[{Gschwind and Drexl(2019)}]{gschwind2019adaptive}
\bibinfo{author}{Gschwind, T.}, \bibinfo{author}{Drexl, M.}, \bibinfo{year}{2019}.
\newblock \bibinfo{title}{Adaptive large neighborhood search with a constant-time feasibility test for the dial-a-ride problem}.
\newblock \bibinfo{journal}{Transportation Science} \bibinfo{volume}{53}, \bibinfo{pages}{480--491}.
\bibitem[{Gschwind and Irnich(2015)}]{gschwind2015effective}
\bibinfo{author}{Gschwind, T.}, \bibinfo{author}{Irnich, S.}, \bibinfo{year}{2015}.
\newblock \bibinfo{title}{Effective handling of dynamic time windows and its application to solving the dial-a-ride problem}.
\newblock \bibinfo{journal}{Transportation Science} \bibinfo{volume}{49}, \bibinfo{pages}{335--354}.
\bibitem[{Ho et~al.(2018)Ho, Szeto, Kuo, Leung, Petering and Tou}]{ho2018survey}
\bibinfo{author}{Ho, S.C.}, \bibinfo{author}{Szeto, W.Y.}, \bibinfo{author}{Kuo, Y.H.}, \bibinfo{author}{Leung, J.M.}, \bibinfo{author}{Petering, M.}, \bibinfo{author}{Tou, T.W.}, \bibinfo{year}{2018}.
\newblock \bibinfo{title}{A survey of dial-a-ride problems: Literature review and recent developments}.
\newblock \bibinfo{journal}{Transportation Research Part B: Methodological} \bibinfo{volume}{111}, \bibinfo{pages}{395--421}.
\bibitem[{Huangfu and Hall(2018)}]{huangfu2018parallelizing}
\bibinfo{author}{Huangfu, Q.}, \bibinfo{author}{Hall, J.J.}, \bibinfo{year}{2018}.
\newblock \bibinfo{title}{Parallelizing the dual revised simplex method}.
\newblock \bibinfo{journal}{Mathematical Programming Computation} \bibinfo{volume}{10}, \bibinfo{pages}{119--142}.
\bibitem[{Jain and Hentenryck(2011)}]{DBLP:conf/cp/JainH11}
\bibinfo{author}{Jain, S.}, \bibinfo{author}{Hentenryck, P.V.}, \bibinfo{year}{2011}.
\newblock \bibinfo{title}{Large neighborhood search for dial-a-ride problems}, in: \bibinfo{editor}{Lee, J.H.} (Ed.), \bibinfo{booktitle}{Principles and Practice of Constraint Programming - {CP} 2011 - 17th International Conference, {CP} 2011, Perugia, Italy, September 12-16, 2011. Proceedings}, \bibinfo{publisher}{Springer}.  \bibinfo{pages}{400--413}.
\newblock \doi{10.1007/978-3-642-23786-7\_31}.
\bibitem[{Jorgensen et~al.(2007)Jorgensen, Larsen and Bergvinsdottir}]{jorgensen2007solving}
\bibinfo{author}{Jorgensen, R.M.}, \bibinfo{author}{Larsen, J.}, \bibinfo{author}{Bergvinsdottir, K.B.}, \bibinfo{year}{2007}.
\newblock \bibinfo{title}{Solving the dial-a-ride problem using genetic algorithms}.
\newblock \bibinfo{journal}{Journal of the operational research society} \bibinfo{volume}{58}, \bibinfo{pages}{1321--1331}.
\bibitem[{Joshi et~al.(2019)Joshi, Laurent and Bresson}]{joshi2019efficient}
\bibinfo{author}{Joshi, C.K.}, \bibinfo{author}{Laurent, T.}, \bibinfo{author}{Bresson, X.}, \bibinfo{year}{2019}.
\newblock \bibinfo{title}{An efficient graph convolutional network technique for the travelling salesman problem}.
\newblock \bibinfo{journal}{arXiv preprint arXiv:1906.01227}.
\bibitem[{Kang et~al.(2019)Kang, Chen and Meng}]{kang2019bus}
\bibinfo{author}{Kang, L.}, \bibinfo{author}{Chen, S.}, \bibinfo{author}{Meng, Q.}, \bibinfo{year}{2019}.
\newblock \bibinfo{title}{Bus and driver scheduling with mealtime windows for a single public bus route}.
\newblock \bibinfo{journal}{Transportation Research Part C: Emerging Technologies} \bibinfo{volume}{101}, \bibinfo{pages}{145--160}.
\bibitem[{Karimi-Mamaghan et~al.(2022)Karimi-Mamaghan, Mohammadi, Meyer, Karimi-Mamaghan and Talbi}]{karimi2022machine}
\bibinfo{author}{Karimi-Mamaghan, M.}, \bibinfo{author}{Mohammadi, M.}, \bibinfo{author}{Meyer, P.}, \bibinfo{author}{Karimi-Mamaghan, A.M.}, \bibinfo{author}{Talbi, E.G.}, \bibinfo{year}{2022}.
\newblock \bibinfo{title}{Machine learning at the service of meta-heuristics for solving combinatorial optimization problems: A state-of-the-art}.
\newblock \bibinfo{journal}{European Journal of Operational Research} \bibinfo{volume}{296}, \bibinfo{pages}{393--422}.
\bibitem[{Kingma and Ba(2014)}]{kingma2014adam}
\bibinfo{author}{Kingma, D.P.}, \bibinfo{author}{Ba, J.}, \bibinfo{year}{2014}.
\newblock \bibinfo{title}{Adam: A method for stochastic optimization}.
\newblock \bibinfo{journal}{arXiv preprint arXiv:1412.6980}.
\bibitem[{Kipf and Welling(2016)}]{kipf2016semi}
\bibinfo{author}{Kipf, T.N.}, \bibinfo{author}{Welling, M.}, \bibinfo{year}{2016}.
\newblock \bibinfo{title}{Semi-supervised classification with graph convolutional networks}.
\newblock \bibinfo{journal}{arXiv preprint arXiv:1609.02907}.
\bibitem[{Kirchler and Calvo(2013)}]{kirchler2013granular}
\bibinfo{author}{Kirchler, D.}, \bibinfo{author}{Calvo, R.W.}, \bibinfo{year}{2013}.
\newblock \bibinfo{title}{A granular tabu search algorithm for the dial-a-ride problem}.
\newblock \bibinfo{journal}{Transportation Research Part B: Methodological} \bibinfo{volume}{56}, \bibinfo{pages}{120--135}.
\bibitem[{Kotary et~al.(2021)Kotary, Fioretto, Van~Hentenryck and Wilder}]{kotary2021end}
\bibinfo{author}{Kotary, J.}, \bibinfo{author}{Fioretto, F.}, \bibinfo{author}{Van~Hentenryck, P.}, \bibinfo{author}{Wilder, B.}, \bibinfo{year}{2021}.
\newblock \bibinfo{title}{End-to-end constrained optimization learning: A survey}.
\newblock \bibinfo{journal}{arXiv preprint arXiv:2103.16378}.
\bibitem[{Lu et~al.(2022)Lu, Nie, Mahmoudi, Ou, Li and Zhou}]{lu2022rich}
\bibinfo{author}{Lu, J.}, \bibinfo{author}{Nie, Q.}, \bibinfo{author}{Mahmoudi, M.}, \bibinfo{author}{Ou, J.}, \bibinfo{author}{Li, C.}, \bibinfo{author}{Zhou, X.S.}, \bibinfo{year}{2022}.
\newblock \bibinfo{title}{Rich arc routing problem in city logistics: Models and solution algorithms using a fluid queue-based time-dependent travel time representation}.
\newblock \bibinfo{journal}{Transportation Research Part B: Methodological} \bibinfo{volume}{166}, \bibinfo{pages}{143--182}.
\bibitem[{Mazyavkina et~al.(2021)Mazyavkina, Sviridov, Ivanov and Burnaev}]{mazyavkina2021reinforcement}
\bibinfo{author}{Mazyavkina, N.}, \bibinfo{author}{Sviridov, S.}, \bibinfo{author}{Ivanov, S.}, \bibinfo{author}{Burnaev, E.}, \bibinfo{year}{2021}.
\newblock \bibinfo{title}{Reinforcement learning for combinatorial optimization: A survey}.
\newblock \bibinfo{journal}{Computers \& Operations Research} \bibinfo{volume}{134}, \bibinfo{pages}{105400}.
\bibitem[{Morabit et~al.(2021)Morabit, Desaulniers and Lodi}]{morabit2021machine}
\bibinfo{author}{Morabit, M.}, \bibinfo{author}{Desaulniers, G.}, \bibinfo{author}{Lodi, A.}, \bibinfo{year}{2021}.
\newblock \bibinfo{title}{Machine-learning--based column selection for column generation}.
\newblock \bibinfo{journal}{Transportation Science} \bibinfo{volume}{55}, \bibinfo{pages}{815--831}.
\bibitem[{Morabit et~al.(2023)Morabit, Desaulniers and Lodi}]{morabit2023machine}
\bibinfo{author}{Morabit, M.}, \bibinfo{author}{Desaulniers, G.}, \bibinfo{author}{Lodi, A.}, \bibinfo{year}{2023}.
\newblock \bibinfo{title}{Machine-learning--based arc selection for constrained shortest path problems in column generation}.
\newblock \bibinfo{journal}{INFORMS Journal on Optimization} \bibinfo{volume}{5}, \bibinfo{pages}{191--210}.
\bibitem[{Nazari et~al.(2018)Nazari, Oroojlooy, Snyder and Tak{\'a}c}]{nazari2018reinforcement}
\bibinfo{author}{Nazari, M.}, \bibinfo{author}{Oroojlooy, A.}, \bibinfo{author}{Snyder, L.}, \bibinfo{author}{Tak{\'a}c, M.}, \bibinfo{year}{2018}.
\newblock \bibinfo{title}{Reinforcement learning for solving the vehicle routing problem}.
\newblock \bibinfo{journal}{Advances in neural information processing systems} \bibinfo{volume}{31}.
\bibitem[{{Partnership for an Advanced Computing Environment (PACE)}(2017)}]{PACE2017}
\bibinfo{author}{{Partnership for an Advanced Computing Environment (PACE)}}, \bibinfo{year}{2017}.
\newblock \bibinfo{title}{{PACE - Partnership for an Advanced Computing Environment}}.
\newblock \bibinfo{note}{[Online; accessed 18-December-2023]}.
\bibitem[{Portugal et~al.(2009)Portugal, Louren{\c{c}}o and Paix{\~a}o}]{portugal2009driver}
\bibinfo{author}{Portugal, R.}, \bibinfo{author}{Louren{\c{c}}o, H.R.}, \bibinfo{author}{Paix{\~a}o, J.P.}, \bibinfo{year}{2009}.
\newblock \bibinfo{title}{Driver scheduling problem modelling}.
\newblock \bibinfo{journal}{Public Transport} \bibinfo{volume}{1}, \bibinfo{pages}{103--120}.
\bibitem[{PyTorch-Geometric(2024)}]{pytorch_geometric_batching}
\bibinfo{author}{PyTorch-Geometric}, \bibinfo{year}{2024}.
\newblock \bibinfo{title}{Batching — pytorch geometric 2023 documentation}.
\newblock \bibinfo{howpublished}{\url{https://pytorch-geometric.readthedocs.io/en/latest/advanced/batching.html}}.
\newblock \bibinfo{note}{Accessed: 2024-12-27}.
\bibitem[{Qu and Bard(2015)}]{qu2015branch}
\bibinfo{author}{Qu, Y.}, \bibinfo{author}{Bard, J.F.}, \bibinfo{year}{2015}.
\newblock \bibinfo{title}{A branch-and-price-and-cut algorithm for heterogeneous pickup and delivery problems with configurable vehicle capacity}.
\newblock \bibinfo{journal}{Transportation Science} \bibinfo{volume}{49}, \bibinfo{pages}{254--270}.
\bibitem[{Riley et~al.(2019)Riley, Legrain and Van~Hentenryck}]{riley2019column}
\bibinfo{author}{Riley, C.}, \bibinfo{author}{Legrain, A.}, \bibinfo{author}{Van~Hentenryck, P.}, \bibinfo{year}{2019}.
\newblock \bibinfo{title}{Column generation for real-time ride-sharing operations}, in: \bibinfo{booktitle}{Integration of Constraint Programming, Artificial Intelligence, and Operations Research: 16th International Conference, CPAIOR 2019, Thessaloniki, Greece, June 4--7, 2019, Proceedings 16}, \bibinfo{organization}{Springer}.  \bibinfo{pages}{472--487}.
\bibitem[{Ropke et~al.(2007)Ropke, Cordeau and Laporte}]{ropke2007models}
\bibinfo{author}{Ropke, S.}, \bibinfo{author}{Cordeau, J.F.}, \bibinfo{author}{Laporte, G.}, \bibinfo{year}{2007}.
\newblock \bibinfo{title}{Models and branch-and-cut algorithms for pickup and delivery problems with time windows}.
\newblock \bibinfo{journal}{Networks: An International Journal} \bibinfo{volume}{49}, \bibinfo{pages}{258--272}.
\bibitem[{Ropke and Pisinger(2006)}]{ropke2006adaptive}
\bibinfo{author}{Ropke, S.}, \bibinfo{author}{Pisinger, D.}, \bibinfo{year}{2006}.
\newblock \bibinfo{title}{An adaptive large neighborhood search heuristic for the pickup and delivery problem with time windows}.
\newblock \bibinfo{journal}{Transportation science} \bibinfo{volume}{40}, \bibinfo{pages}{455--472}.
\bibitem[{Shen et~al.(2022)Shen, Sun, Li, Eberhard and Ernst}]{shen2022enhancing}
\bibinfo{author}{Shen, Y.}, \bibinfo{author}{Sun, Y.}, \bibinfo{author}{Li, X.}, \bibinfo{author}{Eberhard, A.}, \bibinfo{author}{Ernst, A.}, \bibinfo{year}{2022}.
\newblock \bibinfo{title}{Enhancing column generation by a machine-learning-based pricing heuristic for graph coloring}, in: \bibinfo{booktitle}{Proceedings of the AAAI Conference on Artificial Intelligence},  \bibinfo{pages}{9926--9934}.
\bibitem[{T{\'o}th and Kr{\'e}sz(2013)}]{toth2013efficient}
\bibinfo{author}{T{\'o}th, A.}, \bibinfo{author}{Kr{\'e}sz, M.}, \bibinfo{year}{2013}.
\newblock \bibinfo{title}{An efficient solution approach for real-world driver scheduling problems in urban bus transportation}.
\newblock \bibinfo{journal}{Central European Journal of Operations Research} \bibinfo{volume}{21}, \bibinfo{pages}{75--94}.
\bibitem[{Toth and Vigo(2014)}]{toth2014vehicle}
\bibinfo{author}{Toth, P.}, \bibinfo{author}{Vigo, D.}, \bibinfo{year}{2014}.
\newblock \bibinfo{title}{Vehicle routing: problems, methods, and applications}.
\newblock \bibinfo{publisher}{SIAM}.
\bibitem[{Veli{\v{c}}kovi{\'c} et~al.(2017)Veli{\v{c}}kovi{\'c}, Cucurull, Casanova, Romero, Lio and Bengio}]{velivckovic2017graph}
\bibinfo{author}{Veli{\v{c}}kovi{\'c}, P.}, \bibinfo{author}{Cucurull, G.}, \bibinfo{author}{Casanova, A.}, \bibinfo{author}{Romero, A.}, \bibinfo{author}{Lio, P.}, \bibinfo{author}{Bengio, Y.}, \bibinfo{year}{2017}.
\newblock \bibinfo{title}{Graph attention networks}.
\newblock \bibinfo{journal}{arXiv preprint arXiv:1710.10903}.
\bibitem[{Vinyals et~al.(2015)Vinyals, Fortunato and Jaitly}]{vinyals2015pointer}
\bibinfo{author}{Vinyals, O.}, \bibinfo{author}{Fortunato, M.}, \bibinfo{author}{Jaitly, N.}, \bibinfo{year}{2015}.
\newblock \bibinfo{title}{Pointer networks}.
\newblock \bibinfo{journal}{Advances in neural information processing systems} \bibinfo{volume}{28}.
\bibitem[{Wren and Wren(1995)}]{wren1995genetic}
\bibinfo{author}{Wren, A.}, \bibinfo{author}{Wren, D.O.}, \bibinfo{year}{1995}.
\newblock \bibinfo{title}{A genetic algorithm for public transport driver scheduling}.
\newblock \bibinfo{journal}{Computers \& Operations Research} \bibinfo{volume}{22}, \bibinfo{pages}{101--110}.
\bibitem[{Ying et~al.(2019)Ying, Bourgeois, You, Zitnik and Leskovec}]{ying2019gnnexplainer}
\bibinfo{author}{Ying, Z.}, \bibinfo{author}{Bourgeois, D.}, \bibinfo{author}{You, J.}, \bibinfo{author}{Zitnik, M.}, \bibinfo{author}{Leskovec, J.}, \bibinfo{year}{2019}.
\newblock \bibinfo{title}{Gnnexplainer: Generating explanations for graph neural networks}.
\newblock \bibinfo{journal}{Advances in neural information processing systems} \bibinfo{volume}{32}.
\bibitem[{Yuan et~al.(2022a)Yuan, Chen and Hentenryck}]{DBLP:journals/jair/YuanCH22}
\bibinfo{author}{Yuan, E.}, \bibinfo{author}{Chen, W.}, \bibinfo{author}{Hentenryck, P.V.}, \bibinfo{year}{2022}a.
\newblock \bibinfo{title}{Reinforcement learning from optimization proxy for ride-hailing vehicle relocation}.
\newblock \bibinfo{journal}{J. Artif. Intell. Res.} \bibinfo{volume}{75}, \bibinfo{pages}{985--1002}.
\newblock \doi{10.1613/JAIR.1.13794}.
\bibitem[{Yuan et~al.(2022b)Yuan, Jiang and Song}]{yuan2022neural}
\bibinfo{author}{Yuan, H.}, \bibinfo{author}{Jiang, P.}, \bibinfo{author}{Song, S.}, \bibinfo{year}{2022}b.
\newblock \bibinfo{title}{The neural-prediction based acceleration algorithm of column generation for graph-based set covering problems}, in: \bibinfo{booktitle}{2022 IEEE International Conference on Systems, Man, and Cybernetics (SMC)}, \bibinfo{organization}{IEEE}.  \bibinfo{pages}{1115--1120}.

\end{thebibliography}

\end{document}